\documentclass[12pt]{article}
\usepackage{amsfonts,amssymb}
\usepackage{url}
\usepackage{graphicx,color}
\newtheorem{defin}{Definition}[section]

\newtheorem{lem}[defin]{Lemma}
\newtheorem{lemma}[defin]{Lemma}

\newtheorem{prob}[defin]{Problem}
\newtheorem{problem}[defin]{Problem}
\newtheorem{thm}[defin]{Theorem}
\newtheorem{theorem}[defin]{Theorem}
\newtheorem{cor}[defin]{Corollary}

\newtheorem{clm}{Claim}

\newtheorem{con}[defin]{Conjecture}

\def\qed{\hbox{\kern1pt\vrule height6pt width4pt
depth1pt\kern1pt}\medskip}
\def\bproof{\par\noindent{\bf Proof.\enspace}\rm}

\newcommand{\real}{{\mathbb{R}}}
\newcommand{\complex}{{\mathbb{C}}}
\newcommand{\rat}{{\mathbb{Q}}}

\newcommand{\scrrm}{{\mathcal{M}}}
\newcommand{\scrqs}{{\mathcal{QS}}}

\newcommand{\sm}{\setminus}
\newcommand{\tran}{{\mbox{\rm td}}}
\newcommand{\eproof}{\hfill $\bullet$\\}
\newcommand{\rank}{\mbox{\rm rank }}
\newcommand{\rankg}{\mbox{\rm rank}}
\newcommand{\Arg}{\mbox{\rm Arg }}

\title{The number of equivalent realisations of a rigid graph}
\author{
Bill Jackson \thanks{School of Mathematical Sciences, Queen Mary
University of London, Mile End Road, London E1 4NS, England. E-mail:
b.jackson@qmul.ac.uk} \and J.C. Owen \thanks{Siemens, Park House,
Cambridge CB3 0DU, England. E-mail: owen.john.ext@siemens.com} }
\begin{document}
\maketitle

\begin{abstract}
Given a rigid realisation of a graph $G$ in ${\mathbb R}^2$, it is an open problem to determine the maximum number of
pairwise non-congruent realisations which have the same edge lengths as the given realisation.
This problem can be restated as finding the number of solutions of a related system of quadratic equations and in
this context it is natural to consider the number of solutions in ${\mathbb C}^2$ rather that ${\mathbb R}^2$.
We show that the number of complex solutions, $c(G)$, is the same for all generic realisations of a rigid graph $G$,
characterise the graphs $G$ for which $c(G)=1$, and show that the problem of determining $c(G)$ can be reduced to the
case when $G$ is $3$-connected and has no non-trivial $3$-edge-cuts. We consider the effect of the
Henneberg moves and the vertex-splitting operation on $c(G)$. We use our results to determine $c(G)$ exactly for two important families of graphs, and show that the graphs in both families have $c(G)$ pairwise equivalent generic real realisations. We also show that every planar isostatic graph on $n$ vertices has at least $2^{n-3}$ pairwise equivalent real realisations.
\end{abstract}

\section{Introduction}

Graphs with geometrical constraints provide natural models for a
variety of applications, including Computer-Aided Design, sensor
networks and flexibility in molecules. Given a graph $G$ and

prescribed lengths for its edges, a basic problem is to determine
whether $G$ has a straight line realisation in Euclidean
$d$-dimensional space with these given lengths.
Closely related problems are to determine whether a
given realisation is unique or, more generally,  determine how many distinct
realisations exist with the same edge lengths. Saxe \cite{S} has
shown that both the existence and uniqueness problems are NP-hard.
However, this hardness relies on algebraic relations between
coordinates of vertices, and for practical purposes it is natural
to study generic realisations.

Gortler, Healy and Thurston \cite{GHT} showed
that the uniqueness of a generic realisation in $\real^d$ depends only on the
structure of the underlying graph, and we say that a graph $G$ is {\em globally rigid} in $\real^d$ if it has a unique generic realisation in $\real^d$.
It can be seen that $G$ is globally rigid in $\real$ if and only if $G$ is
equal to $K_2$ or is 2-connected. Globally rigid graphs  in $\real^2$ are characterised by a combination of
results due to Hendrickson \cite{H}, Connelly \cite{C}, and Jackson
and Jord\'an \cite{JJ}. No characterisations are known in $\real^d$
when $d\geq 3$.

\begin{figure}
\vspace{-2cm}
\centering
\includegraphics[scale=0.7]{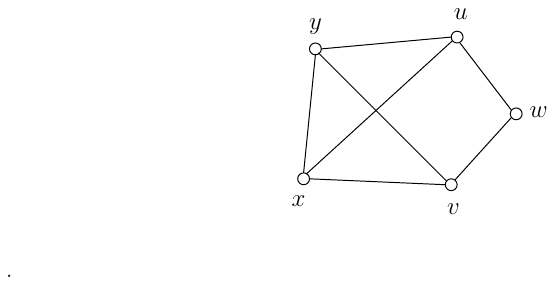}
\vspace{-15cm}
\caption{A realisation of a graph $G$ in $\real^2$. The only other equivalent
realisation is obtained by reflecting the vertex $w$ in the line
through $\{u,v\}$.} \label{fig3}
\end{figure}

\begin{figure}
\vspace{-2cm}
\centering
\includegraphics[scale=0.7]{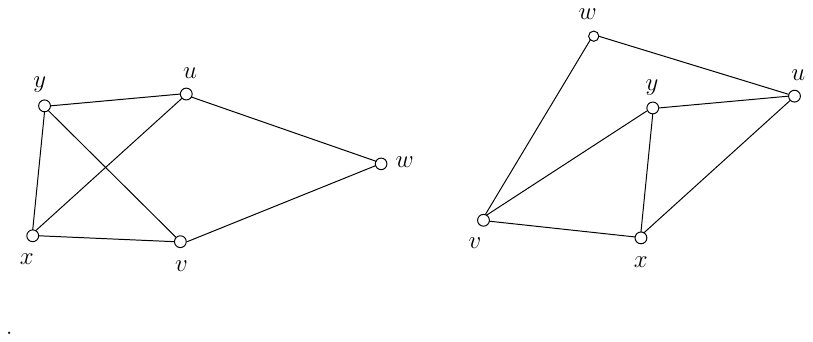}
\vspace{-15cm}
\caption{Two equivalent realisations of the graph $G$ of Figure
\ref{fig3} in $\real^2$. Two other equivalent realisation can be obtained from
these by reflecting the vertex $w$ in the line through $\{u,v\}$,
giving four different equivalent realisations in $\real^2$.} \label{fig4}\end{figure}

In contrast, the number of realisations which are equivalent to,
i.e. have the same edge lengths as, a given generic realisation of a
graph in $\real^d$ may depend on both the graph and the realisation
when $d\geq 2$, see Figures \ref{fig3} and \ref{fig4}. Bounds on the
maximum number of equivalent realisations, where the maximum is
taken over all possible realisations of a given graph,
are obtained by Borcea and Streinu in \cite{BS}, and this number is
determined exactly for generic realisations of an important family of graphs by Jackson,
Jord\'an, and Szabadka in \cite{JJS}.

The set of all realisations which are equivalent to a given
realisation can be represented as the set of solutions to a system
of quadratic equations. In this setting it is natural to consider
the number of complex solutions. This number gives an upper bound on
the number of real solutions which often plays a crucial role in  calculating  the exact number of real solutions, see for example \cite{ETV, EM, ST}.
The number of complex solutions is also much better
behaved than the number of real solutions. For example, we shall show that the number of complex solutions is the same for all
generic realisation of given graph. The realisations of the graph $G$ shown in Figures \ref{fig3} and \ref{fig4} both have four equivalent complex realisations. Only two of these are real in Figure  \ref{fig3}, but  all four are real in Figure \ref{fig4}.

Gortler and Thurston \cite{GT} recently showed that a graph has a unique generic realisation
in $\complex^d$ if and only if it is globally rigid in $\real^d$. This implies that the above mentioned characterisations of
globally rigid graphs in $\real^d$ for $d=1,2$ extend immediately to $\complex^d$, and
explains the apparent inconsistency that having a unique real realisation is a generic
property whereas the number of different real realisations is not.

We will concentrate on the 2-dimensional case in this paper.
We show that the number, $c(G)$, of complex realisations of a rigid
graph $G$ which are equivalent to a given generic realisation is
finite and is the same for all generic realisations. We then
consider the effect of graph operations on $c(G)$. It is known that
a type 1 Henneberg move doubles $c(G)$. We will show that a type 2 Henneberg move on a
redundant edge does not increase $c(G)$, and use this to give a short proof that the characterization of graphs with unique generic realisations in $\real^2$ extends
to $\complex^2$.
We next show that the vertex splitting move increases $c(G)$ by  a factor of at least two when $G$ is isostatic and that the same result holds for the maximum number of pairwise equivalent generic real realisations of $G$. We use this to deduce that every planar isostatic graph on $n$ vertices has a generic real realisation with at least $2^{n-3}$ equivalent realisations. We next consider operations which glue two graphs $G_1,G_2$ together by either associating two pairs of vertices in each graph or by adding three edges between them, and show how
$c(G)$ can be computed from $c(G_1)$ and $c(G_2)$. We use these
results to determine $c(G)$ for a family of quadratically solvable graphs and for graphs with a connected rigidity matroid, and show that the graphs in both families have $c(G)$ pairwise equivalent generic real realisations. We also show that every planar isostatic graph on $n$ vertices has at least $2^{n-3}$ pairwise equivalent real realisations. We
close with a short section of examples and open problems.

\section{Definitions and notation}

A {\em complex (real) realisation} of a graph $G=(V,E)$ is a map $p$
from $V$ to $\complex^2$ ($\real^2$). We also refer to the
ordered pair $(G,p)$ as a {\em framework}. The coordinates  of a point
$p(v)$ are $x(p(v))$ and $y(p(v))$. A framework $(G,p)$ is
{\em generic} if the set of all coordinates of the points $p(v)$,
$v\in V$, is algebraically independent over $\rat$.

For $P=(x,y)\in \complex^2$ let $d(P)=x^2+y^2$ and $\|P\|=(|x|^2+|y|^2)^{1/2}$, where $|.|$ denotes the modulus of a complex number. Two frameworks
$(G,p)$ and $(G,q)$ are {\em equivalent} if
$d(p(u)-p(v))=d(q(u)-q(v))$ for all $uv\in E$, and are {\em
congruent} if $d(p(u)-p(v))=d(q(u)-q(v))$ for all $u,v\in V$.

A framework $(G,p)$ is {\em complex}, respectively {\em real},  {\em rigid} if there exists an
$\epsilon>0$ such that every complex, respectively real, framework $(G,q)$ which is
equivalent to $(G,p)$ and satisfies
$\|(p(v)-q(v)\|<\epsilon$ for all $v\in V$, is
congruent to $(G,p)$.
 Equivalently, every continuous motion of the
points $p(v)$, $v\in V$, in $\complex^2$, respectively $\real^2$, which respects the length
constraints results in a framework which is congruent to $(G,p)$. Note that real rigidity
considers only  the $real$ frameworks
which are equivalent to a given real framework,
whereas complex rigidity
considers all  equivalent complex frameworks,
some or all of which may in fact be real.

The {\it rigidity matrix} of a framework $(G,p)$ is the matrix
$R(G,p)$ of size $|E|\times 2|V|$, where, for each edge $v_iv_j\in
E$, in the row corresponding to $v_iv_j$, the entries in the two
columns corresponding to vertices $v_i$ and $v_j$ contain the two
coordinates of $(p(v_i)-p(v_j))$ and $(p(v_j)-p(v_i))$,
respectively, and the remaining entries are zeros. The framework is {\em infinitessimally
rigid} if  $\rank R(G,p)=2|V|-3$.\footnote{We always have $\rank R(G,p)\leq 2|V|-3$ since its null space
always contains three linearly independent vectors coresponding to two translations and a rotation of the framework.}
    Asimow and Roth
\cite{AR} showed that infinitessimal rigidity is a sufficient
condition for the real rigidity of $(G,p)$, and that the two properties are equivalent
when $(G,p)$ is generic. This implies that real rigidity is a
generic property and we say that $G$ is {\em
rigid} if some/every generic real realisation of $G$ is real rigid. Theorem 5.7 below implies that complex rigidity is also a
generic property and that a graph $G$ is complex rigid  if and only if it is real rigid. This allows us to describe  a graph as
being {\em rigid} without the need to distinguish between real and complex rigidity.
Rigid graphs are characterised by results of Laman \cite{L} and
Lov\'asz and Yemini \cite{LY}. We refer the reader to \cite{Whlong} for more information on the
rigidity of graphs.

Given a  complex or real framework $(G,p)$,
the fact that an algebraic variety can only contain finitely many isolated points implies that the maximum number of pairwise non-congruent rigid frameworks which are equivalent to $(G,p)$ is finite.
We denote the number of such complex, or real, frameworks by  $c(G,p)$, and  $r(G,p)$,  respectively.
We will mostly be concerned with the case when $G$ is rigid and  $(G,p)$ is generic. In this case all equivalent frameworks are rigid and hence
$c(G,p)$, and  $r(G,p)$, will count the total number of non-congruent equivalent complex, and real, frameworks.


\section{Congruent realisations}
Given a complex realisation of a rigid graph it will be useful to
have a `canonical representative' for each congruence class in the
set of all equivalent realisations. The following lemmas will enable
us to do this.

\begin{lemma} \label{basic}(a) Let $P_0\in \complex^2$ and
$\tau:\complex^2\to\complex^2$ by $\tau(P)=P+P_0$. Then
$d(P-Q)=d(\tau(P)-\tau(Q))$ for all $P,Q\in\complex^2$.\\
(b) Let $z_1,z_2\in \complex$ such that $z_1^2+z_2^2=1$ and put
$$M=\left(
\begin{array}{rr}
z_1 & z_2\\
-z_2 & z_1
\end{array}
\right)\,.
$$
Let $\rho:\complex^2\to\complex^2$ by $\rho(P)=MP$. Then
$d(P-Q)=d(\rho(P)-\rho(Q))$ for all $P,Q\in\complex^2$.\\
(c) Put
$$N=\left(
\begin{array}{rr}
1 & 0\\
0 & -1
\end{array}
\right)\,.
$$
Let $\theta:\complex^2\to\complex^2$ by $\theta(P)=NP$. Then
$d(P-Q)=d(\theta(P)-\theta(Q))$ for all $P,Q\in\complex^2$.
\end{lemma}
\bproof (a) is immediate since $\tau(P)-\tau(Q)=P-Q$. To prove (b)
and (c), let $P-Q=(a,b)$. Then
\begin{eqnarray*}
d(\rho(P)-\rho(Q))&=&d(MP-MQ)=d(M(P-Q))=d(z_1a+z_2b,-z_2a+z_1b)\\
&=& (z_1a+z_2b)^2+(-z_2a+z_1b)^2 =a^2+b^2=d(P-Q).
\end{eqnarray*}
     since
$z_1^2+z_2^2=1$.
Similarly \begin{eqnarray*}
d(\theta(P)-\theta(Q))&=&d(NP-NQ)=d(N(P-Q))=d(a,-b)=d(P-Q).
\end{eqnarray*}
\eproof

Three distinct points $P_1,P_2,P_3\in \complex^2$ are {\em collinear} if
$P_2-P_1=z(P_3-P_1)$ for some $z\in \complex$.

\begin{lemma}\label{uniquecon}
Let $P_1,P_2,P_3$ be three distinct points in $\complex^2$
which are not collinear.
Suppose that $M,M'$ are $2\times2$ complex matrices, $t,t'\in
\complex^2$, and that $MP_i+t=M'P_i+t'$ for all $1\leq i\leq 3$.
Then $M=M'$ and $t=t'$.
\end{lemma}
\bproof Since $P_1,P_2,P_3$ are not collinear, $P_2-P_1$ and
$P_3-P_1$ are linearly independent. Furthermore
$(M-M')(P_2-P_1)=0=(M-M')(P_3-P_1)$. Hence $M-M'=0$. Thus $M=M'$ and
$t=t'$. \eproof

\begin{lem}\label{standard} Let $(G,p)$ be a complex realisation of a graph
$G=(V,E)$ with $V=\{v_1,v_2,\ldots,v_n\}$ and $n\geq 3$. Then
$(G,p)$ is congruent to a realisation $(G,q)$ with $q(v_1)=(0,0)$.
Furthermore, if $d(p(v_1)-p(v_2))\neq 0$,
then there exists a unique realisation $(G,q^*)$ which is congruent to $(G,p)$ and satisfies
$q^*(v_1)=(0,0)$, $q^*(v_2)=(0,b_2)$ and $q^*(v_3)=(a_3,b_3)$ for
some $b_2,a_3,b_3\in \complex$ with $b_2\neq 0$, $\Arg b_2\in
(0,\pi]$, and either $a_3=0$ or $\Arg a_3\in
(0,\pi]$.
\end{lem}
\bproof Define $(G,q)$ by putting $q(v_i)=p(v_i)-p(v_1)$ for all $v_i\in
V$. Then $q(v_1)=(0,0)$ and $(G,q)$ is congruent to $(G,p)$ by Lemma
\ref{basic}(a).

Now suppose that $d(p(v_2)-p(v_1))\neq 0$. Then there exists a unique $b_2\in \complex\sm \{0\}$ such that  $d(p(v_2)-p(v_1))=b_2^2$ and $\Arg b_2\in
(0,\pi]$. Let
$q(v_2)=(a,b)$. Then
$$a^2+b^2=d(q(v_2)-q(v_1))=d(p(v_2)-p(v_1))=b_2^2.$$
Put $z_1=b/b_2$ and $z_2=-a/b_2$. Then $z_1^2+z_2^2=1$. We may now
define the matrix $M$ as in Lemma \ref{basic}(b) and define a
realisation $(G,\tilde q)$ by putting $\tilde q(v_i)=Mq(v_i)$ for all
$v_i\in V$. We then have $\tilde q(v_1)=q(v_1)=(0,0)$ and $\tilde
q(v_2)=(0,b_2)$.  Then $(G,\tilde q)$ is congruent to $(G,p)$ by
Lemma \ref{basic}(b). Let $\tilde q(v_3)=(a_3,b_3)$. If $a_3=0$ or
$\Arg a_3\in (0,\pi]$ we put $q^*=\tilde q$; if $\Arg a_3\in
(-\pi,0]$ we put $q^*(v_i)=Nq(v_i)$ for all $v_i\in V$, where $N$ is
the matrix defined in Lemma \ref{basic}(c).
By Lemma \ref{basic}(c), $(G,q^*)$ is congruent to $(G,p)$ and
satisfies the conditions on $q^*$ given in the statement of the
lemma.

It remains to show that $(G,q^*)$ is unique. We have already seen that $b_2$ is uniquely determined by $p$. Choose $d_1,d_2\in
\complex$ such that $d(p(v_1)-p(v_3))=d_1$ and
$d(p(v_2)-p(v_3))=d_2$. Since $(G,p)$ and $(G,q^*)$ are congruent,
we have $a_3^2+b_3^2=d_1$ and $a_3^2+(b_3-b_2)^2=d_2$. These
equations imply that $b_3$ and $a_3^2$ are uniquely determined by
$p$. Since we also have $a_3=0$ or $\Arg a_3\in (0,\pi]$, $q^*(v_3)=(a_3,b_3)$
is uniquely determined by $p$.

By applying a similar argument as in the proceeding paragraph to
$v_i$ for all $4\leq i\leq n$, we have $q^*(v_i)=(\pm a_i,b_i)$ for
some fixed $a_i,b_i\in \complex$ which are uniquely determined by
$p$. Furthermore, the facts that $(G,q^*)$ is congruent to $(G,p)$
and $d(a_3-a_i,b_3-b_i)\neq d(a_3+a_i,b_3-b_i)$ whenever $a_i\neq 0$, imply that
$q^*(v_i)$ is also uniquely determined by $p$. Hence $(G,q^*)$ is
unique.
\eproof

We say that a framework $(G,q)$ with $G=(V,E)$, $V=\{v_1,v_2,\ldots,v_n\}$ and $n\geq 3$ is in {\em canonical position (with respect to $v_1,v_2,v_3$)} if
$q(v_1)=(0,0)$, $q(v_2)=(0,b_2)$ with $b_2\neq 0$ and $\Arg b_2\in (0,\pi]$, and $q(v_3)=(a_3,b_3)$ with either $a_3=0$ or $\Arg a_3\in (0,\pi]$.
Lemma \ref{standard}  immediately implies:

\begin{lemma}\label{qgen1.6}
Let $(G,p)$ be a complex realisation of a rigid
graph $G=(V,E)$ with $V=\{v_1,v_2,\ldots,v_n\}$ and $n\geq 3$, and let
$S$ be the set of all equivalent realisations of $G$. Suppose that $d(p(v_1)-p(v_2))\neq 0$. Then each congruence class in $S$ has a unique representative $(G,p)$ which is in canonical
position with respect to $v_1,v_2,v_3$.
Furthermore, if the vertices in $p(V)$ are not collinear,
then each
congruence class in $S$ has exactly four realisations $(G,q)$ with
$q(v_1)=(0,0)$ and $q(v_2)=(0,b_2)$ for some $b_2\in \complex\sm \{0\}$, and exactly two of these realisations have $\Arg b_2\in (0,\pi]$.
\end{lemma}
\eproof

\section{Field extensions}

In this section we obtain some preliminary results on field
extensions of $\rat$.
We will use these results in the next section
to prove a key lemma: if $(G,p)$ is a generic realisation of a rigid
graph $G$, and $(G,q)$ is an equivalent realisation in canonical
position,
then the two field extensions we obtain by adding either the coordinates of the
points $q(v)$, $v\in V$, or the values $d(p(u)-p(v))$,
$uv\in E$, to $\rat$ have the same algebraic closure.

A point ${ x}\in \complex^n$ is {\em generic} if its
components form an algebraically independent set over $\rat$.
Given a field $K$ we use $K[X_1,X_2,\ldots,X_n]$ to denote the ring of
polynomials in
the indeterminates $X_1,X_2,\ldots,X_n$ with coefficients in $K$ and
$K(X_1,X_2,\ldots,X_n)$
to denote its field of fractions.
Given
a multivariate polynomial function $f:\complex^n\to \complex^m$ we
use $df|_{
x}$ to denote the Jacobean matrix of $f$ evaluated at a point
${ x}\in \complex^n$. We will obtain several results concerning
$\rat(p)$ and $\rat(f(p))$ when $p$ is a generic point in $\complex^n$.
These will be applied to a generic realisation $(G,p)$ by taking
$f(p)$ to be the vector of `squared edge lengths' in $(G,p)$.

\begin{lemma}\label{genpoint1}
Let $f:\complex^n\to \complex^m$ by $f({ x})= (f_1({ x}),f_2({
x}),\ldots,f_m({ x}))$, where $f_i\in \rat[X_1,X_2,\ldots,X_n]$ for
$1\leq i\leq m$. Suppose that ${ p}$ is a generic point in
$\complex^n$ and $\rank df|_{ p}=m$. Then $f({ p})$ is a generic
point in $\complex^m$.
\end{lemma}
\bproof
Relabelling if necessary, we may suppose that the first $m$ columns
of $df|_{{ p}}$ are linearly independent. Let ${
p}=(p_1,p_2,\ldots,p_n)$. Define $h:\complex^m\to \complex^m$ by
$h(x_1,x_2,\ldots,x_m)=f(x_1,x_2,\ldots,x_m,p_{m+1},\ldots,p_n)$.
and let ${ p}'=(p_1,p_2,\ldots,p_m)$. Then $h({ p}')=f({ p})$ and
$\mbox{rank }dh|_{{ p}'}=m$.

Let $h({ p}')=(\beta_1,\beta_2,\ldots,\beta_m)$. Suppose that
$g(\beta_1,\beta_2,\ldots,\beta_m)=0$ for some polynomial $g$ with
integer coefficients. Then $g(f_1({ p}),f_2({
p}),\ldots,f_m({ p}))=0$. Since $ p$ is generic, we have
$g(h({ x}))=0$ for all ${ x}\in \complex^m$. By the inverse
function theorem $h$ maps a sufficiently small open neighbourhood
$U$ of ${ p}'$ bijectively onto $h(U)$. Thus, for each ${
y}\in h(U)$, there exists ${ x}\in U$ such that $h({ x})={
y}$. This implies that $g({ y})=g(h({ x}))=0$ for each ${
y}\in h(U)$. Since $g$ is a polynomial map and $h(U)$ is an open
subset of $\complex^m$, we have $g\equiv 0$. Hence $h({ p}')=f({
p})$ is generic. \eproof

Given a point ${ p}\in \complex^n$ we use $\rat({ p})$ to
denote the field extension of $\rat$ by the coordinates of ${
p}$. Given fields $K\subseteq L\subseteq \complex$
the {\em transcendence degree} of $L$ over $K$, $\tran[L:K]$,
is the cardinality of a largest subset of $L$ which is algebraically
independent over $K$, see \cite[Section 18.1]{Stew}. (It follows
from the Steinitz exchange axiom, see \cite[Lemma 18.4]{Stew}, that
every set of elements of $L$ which is algebraically independent over
$K$ can be extended to a set of $\tran[L:K]$ elements which is
algebraically independent over $K$.)
    We use  $\overline K$ to denote the algebraic
closure of $K$ in $\complex$. Note that  $\tran[\overline K:K]=0$.


\begin{lemma}\label{genpoint2}
Let $f:\complex^n\to \complex^n$ by $f({ x})= (f_1({ x}),f_2({
x}),\ldots,f_n({ x}))$, where $f_i\in \rat[X_1,X_2,\ldots,X_n]$ for
$1\leq i\leq n$. Suppose that
$f({ p})$ is a generic point in $\complex^n$ for some point
$p\in\complex^n$. Then $\overline{\rat(f(p))}= \overline{\rat(p)}$.
\end{lemma}
\bproof Since $f_i$ is a polynomial with rational coefficients,
we have $f_i({ p})\in \rat(p)$ for all $1\leq i\leq n$. Thus
$\rat(f(p))\subseteq \rat(p)$. Since $f({ p})$ is generic,
$\tran[\rat(f(p)):\rat]=n$. Since $\rat(f(p))\subseteq \rat(p)$ and
$p\in\complex^n$ we have $\tran[\rat(p):\rat]=n$. Thus
$\overline{\rat(f(p))}\subseteq \overline{\rat(p)}$ and
$\tran[\overline{\rat(f(p))}:\rat]=n=\tran[\overline{\rat(p)}:\rat]$.
Suppose $\overline{\rat(f(p))}\neq \overline{\rat(p)}$, and choose
$\gamma\in \overline{\rat(p)}- \overline{\rat(f(p))}$. Then $\gamma$
is not algebraic over $\rat(f(p))$ so $S=\{\gamma,f_1({ p}),f_2({
p}),\ldots,f_n({ p})\}$ is algebraically independent over $\rat$.
This contradicts the facts that $S\subseteq \overline{\rat(p)}$ and
$\tran[\overline{\rat(p)}:\rat]=n$. \eproof

\begin{lemma}\label{genpoint3}
Let $f:\complex^n\to \complex^m$ by $f({ x})= (f_1({ x}),f_2({
x}),\ldots,f_m({ x}))$, where $f_i\in \rat[X_1,X_2,\ldots,X_n]$ for
$1\leq i\leq m$. Let $p$ be a generic point in $\complex^n$ and
suppose that $\rank df|_{p}=n$. Let
$W=\{q\in\complex^n\,:\,f(q)=f(p)\}$. Then $W$ is finite and
$\overline{\rat(p)}=\overline{\rat(q)}$ for all $q\in W$.
\end{lemma}
\bproof Reordering the components of $f$ if necessary, we may
suppose that the first $n$ rows of $df|_{p}$ are linearly
independent. Let $g:\complex^n\to \complex^n$ by $g({ x})= (f_1({
x}),f_2({ x}),\ldots,f_n({ x}))$. Then $\rank dg|_{p}=n$ and the set
$W'=\{q\in\complex^n\,:\,g(q)=q(p)\}$ is finite by \cite[Theorem
2.3]{OP}. Since $W\subseteq W'$, $W$ is also finite. Furthermore,
Lemma \ref{genpoint1} implies that $g(p)$ is a generic point in
$\complex^n$. Lemma \ref{genpoint2} and the fact that $g(p)=g(q)$
now give
$\overline{\rat(p)}=\overline{\rat(g(p))}=\overline{\rat(q)}$.
\eproof

\begin{lemma}\label{genpoint4}
Let $f:\complex^n\to \complex^m$ by $f({ x})= (f_1({ x}),f_2({
x}),\ldots,f_m({ x}))$, where $f_i\in \rat[X_1,X_2,\ldots,X_n]$ for
$1\leq i\leq m$. For each $y\in \complex^n$ let
$$W(y)=\{z\in\complex^n\,:\,f(z)=f(y)\}.$$
Suppose that $p$ and $q$ are generic
points in $\complex^n$ and that $\rank df|_{p}=n$. Then $W(p)$ and
$W(q)$ are both finite and $|W(p)|=|W(q)|$.
\end{lemma}
\bproof The fact that $W(p)$ and $W(q)$ are finite follows from
Lemma \ref{genpoint3}. Since $p=(p_1,p_2,\ldots,p_n)$ and
$q=(q_1,q_2,\ldots,q_n)$ are both generic, $\rat(p)$ and $\rat(q)$
are both isomorphic to
$\rat(X_1,X_2,\ldots,X_n)$ and  we may define an isomorphism
$\theta: \rat(p)\to \rat(q)$ by putting $\theta(c)=c$ for all $c\in
\rat$ and $\theta(p_i)=q_i$ for all $1\leq i\leq n$.
We may extend $\theta$ to an isomorphism
$\tilde\theta:\overline{\rat(p)}\to
\overline{\rat(q)}$.\footnote{This follows from the fact that there
is an isomorphism between any two algebraically closed fields of the
same transcendence degree over $\rat$ , which takes a given
transcendence basis for the first to one for the second, see for
example the proof of \cite[Proposition 8.16]{M}.} We may then apply
$\tilde\theta$ to each component of $\overline{\rat(p)}\,^n$ to
obtain an isomorphism $\Theta:\overline{\rat(p)}\,^n\to
\overline{\rat(q)}\,^n$.

   Suppose $z\in W(p)$. Then $f(z)=f(p)$ and
Lemma \ref{genpoint3} gives $\overline{\rat(p)}=\overline{\rat(z)}$.
It follows that each component of $z$ belongs to
$\overline{\rat(p)}$ and hence $z\in \overline{\rat(p)}\,^n$. Thus
$W(p)\subseteq \overline{\rat(p)}\,^n$. In addition we have
\begin{eqnarray*}
f(\Theta(z))&=&\left[f_1(\Theta(z)),\ldots,f_m(\Theta(z))\right]
=\left[\theta(f_1(z)),\ldots,\theta(f_m(z))\right]\\
&=&\left[\theta(f_1(p)),\ldots,\theta(f_m(p))\right]
=\left[f_1(\Theta(p)),\ldots,f_m(\Theta(p))\right]\\
&=&\left[f_1(q),\ldots,f_m(q)\right] =f(q)
\end{eqnarray*}
so
$\Theta(z)\in W(q)$. Since $\Theta$ is a bijection, this implies
that $|W(p)|\leq |W(q)|$. By symmetry we also have $|W(q)|\leq
|W(p)|$ and hence $|W(p)|= |W(q)|$.
   \eproof

\begin{lemma}\label{realise0}
Let $X_1,X_2,\ldots,X_n$ and $D_1,D_2,\ldots, D_t$ be indeterminates
and let $f_i\in K[X_1,X_2,\ldots,X_n,D_1,D_2,\ldots, D_t]$ for all
$1\leq i \leq m$, for some field $K$ with $\rat\subseteq K\subseteq
\complex$. For each $d\in \complex^t$ let $V_d=\{x\in
\complex^n\,:\,f_i(x,d)=0 \mbox{ for all }1\leq i\leq m\}$. Then
$V_d\neq \emptyset$ for some $d\in \complex^t$ with
$\tran[K(d):K]=t$ if and only if $V_d\neq \emptyset$ for all $d\in
\complex^t$ with $\tran[K(d):K]=t$.
\end{lemma}
\bproof Let
$I$ be the ideal of $K(D)[X]$ generated by
$\{f_i(X,D)\,:\,1\leq i\leq m\}$.
For each $d\in \complex^t$ with $\tran[K(d):K]=t$
let $I_d$ be the
ideal of $K(d)[X]$ generated by $\{f_i(X,d)\,:\,1\leq i\leq m\}$.
There is an
isomorphism from $K(D)(X)$ to $K(d)(X)$ which maps $I$ onto $I_d$.
Furthermore, Hilbert's Weak Nullstellensatz, see \cite{CLO}, tells us
that $V_d\neq \emptyset$ if and only if $I_d$ contains a non-zero
element of $K(d)$. We may use the above isomorphism to deduce that
$V_d\neq \emptyset$ if and only if $I$ contains a non-zero element
of $K(D)$. The lemma now follows since the latter condition is
independent of the choice of $d$. \eproof


\section{Generic frameworks}

   Let $G=(V,E)$ be a graph and $(G,p)$ be a complex realisation of
$G$. Let $V=\{v_1,v_2,\ldots,v_n\}$ and $E=\{e_1,e_2,\ldots,e_m\}$.
We view $p$ as a point ${ p}=(p(v_1),p(v_2),\ldots,p(v_n))$ in
$\complex^{2n}$.
The {\it rigidity map} $d_G:\complex^{2n}\to \complex^m$ is given by
$d_G({ p})=(\ell(e_1),\ell(e_2),\ldots,\ell(e_m))$, where
$\ell(e_i)=d(p(u)-p(v))$ when $e_i=uv$. Note that the evaluation of
the Jacobian of the rigidity map at the point ${ p}\in
\complex^{2n}$ is twice the rigidity matrix of the framework
$(G,p)$. When $H$ is a subgraph of $G$, we will simplify notation and write $d_H(p)$ rather that
$d_H(p|_H)$.

A framework $(G,p)$ is  said to be {\em quasi-generic}
if it is congruent to a generic framework.

\begin{lemma}\label{qgen1}
Suppose that $(G,p)$ is a quasi-generic complex realisation of a
graph $G$. If the rows of the rigidity matrix of $G$ are linearly
independent then $d_G(p)$ is generic.
\end{lemma}
\bproof Choose a generic framework $(G,q)$ congruent to $(G,p)$.
Since the rows of the rigidity matrix of $G$ are linearly
independent, $\rank d(d_G)|_{ q}=|E|$. Hence Lemma \ref{genpoint1}
implies that $d_G({ q})$ is generic. The lemma now follows since
$d_G({ p})=d_G({ q})$. \eproof

A graph $G=(V,E)$ is {\em isostatic} if it is rigid and has
$|E|=2|V|-3$. Note that if $(G,p)$ is a generic realisation of an isostatic graph then
its rigidty matrix has linearly independent rows so $d_G(p)$ is generic by
Lemma \ref{qgen1}.

Our next result allows us to use Lemma \ref{standard} to choose a canonical representative for each congruence class
in the set of all realisations which are equivalent to a given generic realisation of a rigid graph.

\begin{lemma}\label{qgen1.5}
Suppose that $(G,p)$ is a generic complex realisation of a rigid
graph $G=(V,E)$ where $V=\{v_1,v_2,\ldots,v_n\}$ and $n\geq 3$.
Then
$d(p(v_1)-p(v_2))\neq 0$.
\end{lemma}
\bproof Suppose $d(p(v_1)-p(v_2))= 0$. By Lemma \ref{standard},
there exists a realisation $(G,q)$ of $G$ which is congruent
to $(G,p)$ and has $q(v_1)=(0,0)$. Since $d(
q(v_1)-q(v_2))= d(p(v_1)-p(v_2))=0$ we have $
q(v_2)=(b,\pm ib)$ for some $b\in \complex$. Let $H$ be a spanning
isostatic subgraph of $G$. Lemma \ref{qgen1} implies that $d_H(p)$
is generic. Since $d_H(p)=d_H(q)$,
$\tran[\rat(d_H(q)):\rat]=|E(H)|=2n-3$. Since $\rat(d_H(
q))\subseteq \rat(q)$ we have $\tran[\rat(q):\rat]\geq
2n-3$. Since $i=\sqrt{-1}$ is algebraic over $\rat$, this implies
that the set of coordinates of the points $ q(v_i)$, $3\leq
i\leq n$, is algebraically independent over $\rat$. In particular
$q(v_3)\neq (x,\pm ix)$ for all $x\in \complex$. Lemma
\ref{standard}, now gives us a realisation $(G,q')$ of $G$ which is
congruent to $(G,q)$, has $q'(v_1)=(0,0)$ and
$q'(v_3)=(0,b_3)$ for some $b_3\in \complex$. Furthermore
$$d(q'(v_2))=d(q'(v_2)-q'(v_1))=d(q(v_2)-q(v_1))=d(q(v_2))=0$$
so $q'(v_2)=(c,\pm ic)$ for some $c\in \complex$. This implies that
$\tran[\rat(q'):\rat]\leq 2n-4$, and contradicts the facts that
$d_H(q)=d_H(q')$ and $\rat(d_H(q'))\subset \rat(q')$ so
$\tran[\rat(q'):\rat]\geq \tran[\rat(d_H(q)):\rat]=2n-3$. Hence
$d(q(v_1)-q(v_2))\neq 0$.
\eproof

Our next two results show that if $(G,q)$ is equivalent to a generic
realisation $(G,p)$ of a rigid graph and is in canonical position
then the algebraic closures of $\rat(q)$, $\rat(p)$ and $\rat(d_G(p))$ are the
same.

\begin{lemma}\label{qgen2}
Let $(G,p)$ be a complex realisation of an isostatic graph $G=(V,E)$
with $V=\{v_1,v_2,\ldots,v_n\}$. Suppose that $p(v_1)=(0,0)$,
$p(v_2)=(0,y_2)$, $p(v_i)=(x_i,y_i)$ for $3\leq i\leq n$, and
$d_G({p})$ is generic.
Then  $\overline{\rat({ p})}=\overline{\rat({d_G({ p}}))}$.
\end{lemma}
\bproof  Let $f:\complex^{2n-3}\to \complex^{2n-3}$ be defined by
putting $f(z_1,z_2,\ldots,z_{2n-3})$ equal to
$d_G(0,0,0,z_1,z_2,\ldots,z_{2n-3}).$ Let ${
p}'=(y_2,x_3,y_3,\ldots,x_{n},y_n)$. Then $f({ p}')=d_G({ p})$ is
generic, $\rat({ p})=\rat({ p}')$ and ${\rat({d_G({ p}})})=\rat(f({
p}'))$. Lemma \ref{genpoint2} now implies that
$\overline{\rat({p})}=\overline{\rat({
p'})}=\overline{\rat({f({ p'}}))}=\overline{\rat({d_G({ p}}))}$.
\eproof

\begin{lemma}\label{new}
Let $(G,p)$ be a quasi-generic complex realisation of a rigid graph
$G=(V,E)$ with $V=\{v_1,v_2,\ldots,v_n\}$ and $n\geq 3$. Suppose
that $p(v_1)=(0,0)$ and  $p(v_2)=(0,y)$ for some $y\in \complex$.
Let $(G,q)$ be another realisation of $G$ which is equivalent to
$(G,p)$ and has $q(v_1)=(0,0)$ and  $q(v_2)=(0,z)$ for some $z\in
\complex$. Then
$\overline{\rat(p)}=\overline{\rat(q)}=\overline{\rat({d_G({ p}}))}$
and $\tran[\rat(q):\rat]=2n-3$.
\end{lemma}
\bproof  Choose a spanning isostatic subgraph $H$ of $G$.
Lemma \ref{qgen1} implies that $d_H(p)=d_H(q)$ is
generic. Lemma \ref{qgen2} now gives
$\overline{\rat(p)}=\overline{\rat(d_H(p))}=\overline{\rat(q)}$ and
$\tran[\rat(q):\rat]=\tran[\rat(d_H(p):\rat]=2n-3$.
\eproof

Lemma \ref{new} implies that $\tran[\rat(d_G(p):\rat]=2|V|-3$ for
any generic realisation $(G,p)$ of a rigid graph $G$. Our next
result extends this to all graphs. Given a graph $G$ we use
$\rankg(G)$ to denote the rank of the rigidity matrix of a generic
realisation of $G$. A {\em rigid component} of $G$ is a maximal
rigid subgraph of $G$. It is known that the edge-sets of the rigid
components $H_1,H_2,\ldots,H_t$ of $G$ partition $E(G)$ and that
$\rankg(G)=\sum_{i=1}^t(2|V(H_i)|-3)$, see for example \cite{JJ}.

\begin{lemma}\label{grank}
Let $(G,p)$ be a quasi-generic complex realisation of a graph $G$.
Then $\tran[\rat(d_G(p)):\rat]=\rankg(G)$.
\end{lemma}
\bproof Let $H_1,H_2,\ldots,H_t$ be the rigid  components of $G$. By
Lemma \ref{new},
$\tran[\rat(d_{H_i}(p)):\rat]=2|V(H_i)|-3$ for all $1\leq
i\leq t$. Thus
$$\tran[\rat(d_G(p)):\rat]\leq \sum_{i=1}^t \tran[\rat(d_{H_i}(p)):\rat]
=\sum_{i=1}^t (2|V(H_i)|-3)=\rankg(G).$$ On the other hand, we may
apply Lemma \ref{qgen1} to a spanning subgraph $F$ of $G$ whose edge
set corresponds to a maximal set of linearly independent rows of the
rigidity matrix of $(G,p)$ to deduce that
$\tran[\rat(d_G(p)):\rat]\geq \tran[\rat(d_F(p)):\rat]= \rankg(G)$.
Thus $\tran[\rat(d_G(p)):\rat]=\rankg(G)$. \eproof

We can now show that  the number of pairwise
non-congruent realisations of a rigid graph $G$ which are equivalent
to a given generic realisation is the same for all generic
realisations.

\begin{theorem}\label{c(G)} Suppose $(G,p)$ is a generic complex realisation of a
rigid graph $G=(V,E)$. Let $S$ be the set of all equivalent
realisations of $G$. Then the number of congruence classes in $S$ is
finite. Furthermore, this number is the same for all generic
realisations of $G$.
\end{theorem}
\bproof Let $V=\{v_1,v_2,\ldots,v_n\}$ and let $(G,q)$ be another
generic realisation of $G$.  Let $(G,p^*)$ and $(G,q^*)$ be
realisations in canonical position which are congruent to $(G,p)$
and $(G,q)$ respectively. Let $p^*(v_i)=(x_i,y_i)$ and
$q^*(v_i)=(x_i',y_i')$ for $1\leq i\leq n$. Let
$\tilde{p}=(y_2,x_3,y_3,\ldots,x_{n},y_n)$ and
$\tilde{q}=(y'_2,x'_3,y'_3,\ldots,x'_{n},y'_n)$. Then $\tilde{p}$
and $\tilde{q}$ are generic by Lemma \ref{new}.

Let $f:\complex^{2n-3}\to \complex^{2n-3}$ be defined by putting
$f(z_1,z_2,\ldots,z_{2n-3})$ equal to
$d_G(0,0,0,z_1,z_2,\ldots,z_{2n-3})$. Then $\rank
df|_{\tilde{p}}=\rank d(d_G)|_{p}=2n-3$ since $(G,p)$ is infinitesimally rigid (and hence the only vector in the
null space of $R(G,p^*)$ which has a zero in its first three components is the zero vector).
For each $y\in \complex^n$ let
$W(y)=\{z\in\complex^n\,:\,f(z)=f(y)\}.$ By Lemma \ref{genpoint3},
$W(\tilde{p})$ is finite. Since $W(\tilde{p})=4c(G,p)$ by Lemma
\ref{qgen1.6}, $c(G,p)$ is finite. Since
$|W(\tilde{q})|=|W(\tilde{q})|$ by Lemma \ref{genpoint4}, $c(G,q)$
is also finite and $c(G,p)=c(G,q)$. \eproof

As mentioned in the Introduction, we denote the common value of $c(G,p)$
over all generic realisations of $G$ by $c(G)$.
We close this section by using Lemma \ref{qgen1.6} to obtain a lower bound on $c(G)$
 using a non-generic realisation of $G$. Our proof uses the concept of the {\em multiplicity of an isolated solution of a system of polynomial equations}. We refer the reader to \cite[page 224]{SW} for a formal definition but note that an isolated solution $p\in \complex^n$ of a system of $n$ equations in $n$ variables has multiplicity one if the Jacobean of the system has rank $n$ at $p$ and has multiplicity at least  two if the Jacobean has rank less than $n$.

Let $S$ be the set of all rigid frameworks which are equivalent to a given framework $(G,p)$, $\Omega(G,p)$ be the partition of $S$ into congruence classes and $\Omega'(G,p)$ be the set of all congruence classes in $\Omega(G,p)$ which contain frameworks which are rigid but are not infinitesimally rigid and not collinear. By definition we have $c(G,p)=|\Omega|$. Let $c'(G,p)=|\Omega'|$.

\begin{theorem}\label{iso_gen_lem} Suppose $(G,p)$ is a realisation of an
isostatic graph $G=(V,E)$ with $d(p(v_1)-p(v_2))\neq 0$ for some $v_1v_2\in E$. Then $$c(G) \geq c(G,p)+c'(G,p).$$
\end{theorem}
\bproof
Let $(G,q)$ be a generic realisation of $G$.
Since $G$ is isostatic and $(G,q)$ is generic,  $d_G(q)$  is generic over $\rat$ by Lemma \ref{qgen1}.  Let $W_p$ and $W_q$ be the set of all
$t\in \complex^{2|V|}$ such that $t(v_1)=(0,0)$, the first component of $t(v_2)$ is zero, and $(G,t)$ is equivalent to $(G,p)$, respectively $(G,q)$.
Then Lemma \ref{qgen1.6} implies that $W_q$ is a complex algebraic variety defined for a generic set of parameters $d_G(q)$ and has exactly $4\,c(G)$  points. Let $m_i(q)$ be the number of isolated points of $W_q$ with multiplicity $i$. Then $m_i(q)=0$ for $i\geq 2$ because every framework equivalent to $(G,q)$ is infinitesimally rigid. Hence $\Sigma_{i\geq 1} im_i(q) = 4c(G)$.

Similarly let
$m_i(p)$ be the number of isolated points of $W_p$ with multiplicity $i$. Since  $W_q$ is defined by a set of $|E|$ polynomials
in $2|V|-3$ variables with $|E|=2|V|-3$ and $d_G(p)$ is a specialisation  of $d_G(q)$, we have $\Sigma_{i\geq 1} im_i(p) \leq \Sigma_{i\geq 1} im_i(q)$  by \cite[Theorem 7.1.6]{SW}.

Let $(G,\tilde p)$ be a rigid framework which is equivalent to $(G,p)$. Then a similar argument to that used in the derivation of Lemma \ref{qgen1.6} implies that $W_p$ contains at least two isolated points which are congruent to $\tilde p$ if $\tilde p(V)$ is collinear and  at least four isolated points which are congruent to $\tilde p$ if $\tilde p(V)$ is not collinear.
Furthermore, each isolated point of $W_p$ corresponding to a rigid framework which is not infinitesimally rigid (and in particular each isolated point corresponding to a rigid collinear framework) has multiplicity at least two,
and each isolated point corresponding to an infinitesimally rigid framework has multiplicity one.
Hence
$$4c(G,p)+4c'(G,p)\leq \sum_{i\geq 1} im_i(p) \leq \sum_{i\geq 1} im_i(q) = 4c(G).$$

\eproof
Note that Theorem \ref{iso_gen_lem} also holds when $d(p(v_1)-p(v_2))= 0$ for all $v_1v_2\in E$ since in this case we have $c(G,p)=1$, $c'(G,p)=0$ and $c(G)\geq 1$.

The result \cite[Theorem 7.1.6]{SW} we used in the proof of Theorem \ref{iso_gen_lem} is obtained using homotopic continuation.
A purely algebraic proof
for the case when all frameworks equivalent to $(G,p)$ are rigid can be obtained using \cite[Chapter XI]{HP}.

It is not difficult to construct frameworks which show that strict
inequality can hold in Theorem \ref{iso_gen_lem}. For example label
the vertices of $K_4$ as $v_1,v_2,v_3,v_4$, let $H=K_4-v_3v_4$ and
let $G$ be obtained by adding a new vertex $v_5$ and two new edges
$v_5v_3, v_5v_4$ to $H$. It is straightforward to show that $c(G)=4$, for example by using Lemma
\ref{h1} below. However the realisation $(G,p)$ given by  $p(v_1)=(0,0)$,
$p(v_2)=(0,1)$, $p(v_3)=(1,1)$, $p(v_4)=(-1,1)$, and $p(v_5)=(2,3)$
has $c(G,p)=2$. This follows because every realisation $(H,q)$ which
is equivalent but not congruent to $(H,p|_H)$ has $q(v_3)=q(v_4)$
and hence cannot be extended to a realisation of $G$ which is
equivalent to $(G,p)$ (because $d(p(v_5)-p(v_3)) \neq
d(p(v_5)-p(v_4))$). Thus all realisations equivalent to $(G,p)$ are
extensions of $(H,p|_H)$ and there are exactly two ways to do this.

Note also that the conclusion of Theorem \ref{iso_gen_lem} does not hold for rigid graphs which are not isostatic.
For example, label the vertices of $K_5$ as $v_1,v_2,v_3,v_4,v_5$,
and let $G=K_5-v_4v_5$. Then $c(G)=1$ because $G$ is globally rigid.
On the other hand, any rigid realisation $(G,p)$  with $p(v_1)$, $p(v_2)$
and $p(v_3)$ collinear has $c(G,p)\geq 2$ since we may obtain an
equivalent but non-congruent realisation by reflecting $p(v_4)$ in
the line joining $p(v_1)$, $p(v_2)$ and $p(v_3)$.

\section{Graph construction moves}
We first consider the effect of Henneberg moves on the number of
equivalent complex realisations of a rigid graph. The {\em type 1 Henneberg move} on a graph $H$ adds a new vertex
$v$ and two new edges $vx,vy$ from $v$ to distinct vertices $x,y$ of $H$. The {\em type 2 Henneberg move} deletes an edge
$xy$ from $H$ and adds a new vertex
$v$ and three new edges $vx,vy,vz$ from $v$ to $x,y$ and another vertex $z$ of $H$ distinct from $x,y$.

It is straightforward to show that applying the type 1 move will double the number of realisations, see for example \cite{BS,ST}.

\begin{lemma}\label{h1}   Let $G=(V,E)$ be a rigid graph with at least four vertices,
$v_n\in V$ with $N(v_n)=\{v_1,v_2\}$, and  $H=G-v_n$. Then
$c(G)=2c(H)$.
\end{lemma}

We next consider type 2 moves. We need the following result which is
an extension of \cite[Lemma 4.1]{JJS} to complex frameworks. Its
proof uses ideas from simplified versions of the proof of
\cite[Lemma 4.1]{JJS} given in \cite{N,Sz}.

\begin{lemma}\label{con3}   Let $(G,p)$ be a quasi-generic complex framework
and $v_n\in V$ with $N(v_n)=\{v_1,v_2,v_3\}$. Suppose that $(G,q)$
is a complex realisation of $G$ which is equivalent to $(G,p)$. If
$G-v_n$ is rigid then $d(p(v_i)-p(v_j))=d(q(v_i)-q(v_j))$ for all
$1\leq i<j\leq 3$.
\end{lemma}

\bproof By symmetry we need only show that
$d(p(v_1)-p(v_2))=d(q(v_1)-q(v_2))$. Label the vertices of $G$ as
$v_1,\ldots,v_n$ and put $p(v_i)=p_i=(p_{i,1},p_{i,2})$ and
$q(v_i)=q_i=(q_{i,1},q_{i,2})$ for all $1\leq i\leq n$. Since
$G-v_n$ is rigid and $d(v_n)=3$, $G$ is rigid.
By applying Lemma \ref{qgen1.6} to both $(G,p)$ and $(G,q)$,
we may suppose that $p_{1,1}=p_{1,2}=p_{2,2}=0$ and $q_{1,1}=q_{1,2}=q_{2,2}=0$.
\footnote{We have switched the order of the coordinate axes from that given in Lemma \ref{qgen1.6}
since it makes the remainder of the proof more straightforward.}
Then
$$d(p_1-p_2)-d(q_1-q_2)=p_{2,1}^2-q_{2,1}^2$$
so it will suffice to show that $p_{2,1}^2-q_{2,1}^2=0$.

Let $p'=p|_{H}$, $q'=q|_{H}$, $K=\rat(p')$ and $L=\rat(q')$. Consider the equivalent
frameworks $(G-v_n,p')$ and $(G-v_n,q')$. Applying Lemma \ref{new}
to  $G-v_n$, we have $\overline K=\overline L$.
Thus $q_{2,1},q_{3,1},q_{3,2}\in \overline K$. Since
$(G,q)$ is equivalent to $(G,p)$, we have the following equations.
\begin{eqnarray}
q_{n,1}^2+q_{n,2}^2&=&p_{n,1}^2+p_{n,2}^2       \label{eq:A}\\
(q_{n,1}-q_{2,1})^2+ q_{n,2}^2&=&(p_{n,1}-p_{2,1})^2+ p_{n,2}^2      \label{eq:B}\\
(q_{n,1}-q_{3,1})^2+( q_{n,2}-q_{3,2})^2&=&(p_{n,1}-p_{3,1})^2+(
p_{n,2}-p_{3,2})^2 \label{eq:C}
\end{eqnarray}
Subtracting (\ref{eq:A}) from (\ref{eq:B}) and (\ref{eq:C}) we
obtain
\begin{eqnarray}
q_{n,1}&=&\frac{p_{2,1}}{q_{2,1}}p_{n,1}+\frac{q_{2,1}^2-p_{2,1}^2}{2q_{2,1}} \label{eq:D}\\
q_{n,2}&=&\frac{p_{3,1}}{q_{3,2}}p_{n,1}+\frac{p_{3,2}}{q_{3,2}}p_{n,2}-
\frac{q_{3,1}}{q_{3,2}}q_{n,1}+\frac{q_{3,1}^2-p_{3,1}^2+q_{3,2}^2-p_{3,2}^2}{2q_{3,2}}
\label{eq:E}
\end{eqnarray}
We may use (\ref{eq:D}) to eliminate $q_{n,1}$ from the right hand
side of (\ref{eq:E}) to obtain a matrix equation for $q_n$ of the
form
\begin{equation}\label{egy}
q_n=Ap_n+b
\end{equation}
where $A$ is a $2\times 2$ lower triangular matrix with entries in
$\overline K$ and $b\in {\overline K}\,^2$. Rewriting (\ref{eq:A}) as
$q_n^Tq_n=p_n^Tp_n$ and then substituting for $q_n$ using
(\ref{egy}) we obtain
\begin{equation}\label{ketto}
p_n^T(A^TA-I)p_n+2b^TAp_n+b^Tb=0.
\end{equation}
This is a polynomial equation for the components of $p_n$ with
coefficients in $\overline K$. Since $\tran[\rat(p):\rat]=2n-3$ by Lemma
\ref{new},
     $\{p_{n,1},p_{n,2}\}$ is
algebraically independent over $\overline K$. This implies that the polynomial on the left hand side of (\ref{ketto}) is identically
zero. In particular $A^TA=I$ and, since $A$ is lower triangular, $A$
must be a diagonal matrix with $\pm 1$ entries on the diagonal. In
particular $a_{1,1}=p_{2,1}/q_{2,1}=\pm 1$ and hence
$p_{2,1}^2-q_{2,1}^2=0$. \eproof

\begin{lemma}\label{triangle}   Let $G=(V,E)$ be a rigid graph,
$v_n\in V$ with $N(v_n)=\{v_1,v_2,v_3\}$, and $H=(G-v_n)\cup
\{e_1,e_2,e_3\}$ where $e_1=v_1v_2$, $e_2=v_2v_3$ and $e_3=v_1v_3$.
Suppose that $G-v_n$ is rigid. Then $c(G)=c(H)$.
\end{lemma}
\bproof Let $(G,p)$ be a generic realisation of $G$ and $(G,p')$ be
a realisation which is congruent to $(G,p)$ and in canonical
position. Let $S$ be the set of all realisations $(G,q)$ which are
equivalent to $(G,p)$ and in canonical position. Similarly, let
$S^*$ be the set of all realisations $(H,q^*)$ which are equivalent
to $(H,p|_{H})$ and in canonical position. By Lemma
\ref{qgen1.6}, $|S|=c(G)$ and $|S^*|=c(H)$.

Let $F$ be a complete graph with vertex set $\{v_1,v_2,v_3,v_n\}$.
Then Lemma \ref{con3} implies that $(F,q|_{V(F)})$ is congruent to
$(F,p'|_{V(F)})$ for all $(G,q)\in S$. Lemma \ref{qgen1.6} now gives
$q(v_i)=p'(v_i)$ for all $i\in \{1,2,3,n\}$ and all $(G,q)\in S$. We
may use a similar argument to deduce that $q^*(v_i)=p'(v_i)$ for all
$i\in \{1,2,3\}$ and all $(H,q^*)\in S^*$. This implies that the map
$\theta:S\to S^*$ defined by $\theta(G,q)=(H,q|_{V-v_n})$ for all
$(G,q)\in S$ is a bijection.   Hence $c(G)=|S|=|S^*|=c(H)$. \eproof

\begin{cor}\label{h2}   Let $G=(V,E)$ be a rigid graph,
$v_n\in V$ with $N(v_n)=\{v_1,v_2,v_3\}$, and $H=(G-v_n)+e_1$ where
$e_1=v_1v_2$. Suppose that $G-v_n$ is rigid. Then $c(G)\leq c(H)$.
\end{cor}
\bproof By Lemma \ref{triangle}, $c(G)=c(H\cup \{e_2,e_3\})\leq
c(H)$, where $e_2=v_2v_3$ and $e_3=v_1v_3$. \eproof

An edge $e$ in a rigid graph $G$ is {\em redundant} if $G-e$ is rigid. Corollary \ref{h2} tells us that if we extend a rigid graph $H$ by
performing a Henneberg type 2 move on a redundant edge of $H$ then
we do not increase $c(H)$. On the other hand it is not difficult to construct examples with $c(G)=1$ and $c(H)$ arbitrarily large.

It is an open problem to determine the effect that performing a Henneberg
type 2 move on a non-redundant edge has on $c(H)$.

\begin{problem}\label{h2con}   Do there exist universal constants
$k_1,k_2>0$ such that if
$H$ is a rigid and $G$ is obtained by performing a Henneberg type 2 move on a non-redundant edge of $H$,  then $k_1\,c(H)\leq c(G)\leq k_2\,c(H)$?
\end{problem}

We next consider the operation of {\em vertex splitting} introduced by Whiteley in \cite{Wsplit}. Given a vertex $v$ in a graph $H$, this move constructs a new graph $G$ from $H-v$ by partitioning the neighbours  of $v$ into two sets $N_1,N_2$, adding  two new vertices $v_1,v_2$ joined to $N_1,N_2$ respectively, then adding the edge $v_1v_2$ and another edge $v_2x$ for any $x\in N_1$. Whiteley showed that this move preserves rigidity. Since it also preserves the edge count, $G$ will be isostatic whenever $H$ is isostatic.

\begin{lem} \label{vsplit}
Let $H$ be an isostatic graph on at least three vertices and $G$ be obtained from $H$ by applying the vertex splitting move. Then  $c(G) \geq 2c(H)$.
\end{lem}
\bproof
Suppose that $G$ is obtained from $H$ by splitting $v$ into $v_1,v_2$.
Let $(H,p)$ be a generic realisation of $H$. Construct a realisation $(G,q)$ of $G$ by putting $q(v_1)=q(v_2)=p(v)$ and $q(u)=p(u)$ for all other vertices $u$.
We will obtain a bound on $c(G)$ by applying Lemma \ref{iso_gen_lem} to $(G,q)$. It is straightforward to show that $c(G,q)=c(H,p)$. In addition, no framework which is equivalent to $(G,q)$ is infinitesimally rigid (because the row indexed by $v_1v_2$ in its rigidity matrix is zero) or collinear (because $(H,p)$ is generic). Hence $c'(G,q)=c(G,q)$. Lemma \ref{iso_gen_lem} now gives $$c(G) \geq 2c(G,q)=2c(H,p)=2c(H).$$
\eproof

We can also obtain a lower bound on the number of real generic realisations produced by the vertex split operation. We need the following lemma.

\begin{lem} \label{equiv_set} Suppose that $G=(V,E)$ is isostatic and $(G,s_i)$, $1 \leq i \leq m$, are distinct pairwise equivalent infinitesimally rigid  real frameworks in canonical position with respect to three given vertices $v_1,v_2,v_3$. Then, for each $\epsilon>0$, there exists  distinct quasi-generic  pairwise equivalent real frameworks $(G,t_i)$, $1 \leq i \leq m$, in canonical position with respect to $v_1,v_2,v_3$, which  satisfy $\|t_i-s_i\| < \epsilon$ for all $1 \leq i \leq m$.
\end{lem}
\bproof
Let $V=\{v_1,v_2,\ldots,v_n\}$ and $E=\{e_1,e_2,\ldots,e_{2n-3}\}$.
We associate each vector $\hat p=(b_2,a_3,b_3,\ldots,a_n,b_n)\in \real^{2n-3}$ with a real framework $(G,p)$ where $p(v_1)=(0,0)$,  $p(v_2)=(0,b_2)$ and
$p(v_i)=(a_i,b_i)$ for $i\geq 3$.
We can now define a differentiable map $f_G:\real^{2n-3}\to \real^{2n-3}$ by taking $f_G(\hat p)$ to be the ordered vector of squared edge lengths in the framework $(G,p)$.
The rank of the Jacobean matrix $df_G|_{\hat p}$ is equal to the rank of the rigidity matrix of $(G,p)$ and hence $\rank df_G|_{\hat s_i}=2n-3$ for all $1\leq i\leq m$. The inverse function theorem now implies that we can choose open neighbourhoods $N_i$ of $\hat s_i$ and $N$ of $f_G(\hat s_i)$ in $\real^{2n-3}$ such that $f_G$ maps $N_i$ diffeomorphically onto $N$ for all $1\leq i\leq m$. This allows us to choose a generic point $d\in N$ and points $\hat t_i\in N_i$ such that $\|t_i-s_i\|<\epsilon$ and $f_G(\hat t_i)=d$ for all $1\leq i\leq m$. Since $d$ is generic, Lemma \ref{qgen2} implies that each framework $(G,t_i)$ will be quasi-generic. We can ensure that $(G,t_i)$  is in canonical position with respect to $v_1,v_2,v_3$ by choosing $d$ such that $t_i$ is sufficiently close to $s_i$.
\eproof

\begin{thm} \label{rvsplit}
Let $H$ be an isostatic graph on at least three vertices and let $G$ be obtained from $H$ by applying the vertex splitting move.
Let $(H,p)$ be a generic real realisation of $G$. Then there exists a generic realisation $(G,t)$ of $G$ such that $r(G,t)\geq 2r(H,p)$.
\end{thm}
\bproof
Suppose the vertex splitting move splits $v$ into $v_1,v_2$ and adds edges $v_1v_2$ and $v_1x$.
We may suppose that $(H,p)$ is a quasi-generic real realisation of $H$ in canonical position with respect to $v,x,y$ for some vertex $y$ of $H$. We have $p(v)=(0,0)$ and $p(x)=(0,a)$ with $a$ generic. Note that all frameworks equivalent to $(H,p)$ are quasi-generic and hence infinitesimally rigid by Lemma \ref{new}. Construct a real realisation $(G,q)$ of $G$ by putting $q(v_1)=q(v_2)=p(v)$ and $q|_H=p$, and let $S$ be the set of all equivalent realisations $(G,q_i)$, in canonical position with respect to $v_1,x,y$. It is straightforward to show that $r(G,q)=|S|=r(H,p)$.



Let $\hat G$ be the graph obtained from $G$ by performing a Henneberg type 2 move which deletes the edge $v_1v_2$ and adds a new vertex $w$ and new edges $wv_1, wv_2,wx$. Let $(\hat G,\hat q)$ be the framework obtained by putting $\hat q(w)=(a,0)$ and $\hat q|_G=q$, and let $\hat S$ be the set of all equivalent real realisations $(\hat G,\hat q_i)$ which are in canonical position with respect to $v_1,x,y$ and satisfy $\hat q_i(v_1)=\hat q_i(v_2)=(0,0)$.
Then $(\hat G,\hat q_i)\in \hat S$ if and only if $(G,\hat q_i|_G)\in S$ and $\hat q_i(w)=(\pm a,0)$, so  $|\hat S|=2|S|=2r(H,p)$. In addition each $(\hat G,\hat q_i)\in \hat S$ is infinitesimally rigid. To see this suppose that $m$ is an infinitesimal motion of $(\hat G,\hat q_i)$ with $m(v_1)=m(x)=(0,0)$. Then $m(w)=(0,0)$ and so $m(v_2)=(0,0)$. It follows that $m$ induces an infinitesimal motion of the framework $(H,p_i)$ given by $p_i(v)=(0,0)$ and $p_i|_{H-v}=\hat q_i|_{H-v}$, which is zero on $v$ and $x$. The facts that $(H,p_i)$ is equivalent to $(H,p)$ and that all frameworks equivalent to $(H,p)$ are infinitesimally rigid now tells us that $m$ is identically zero.

 We can now use Lemma \ref{equiv_set} to deduce that, for all $\epsilon>0$, there exists  a set $\hat T$ of infinitesimally rigid, pairwise equivalent, quasi-generic real frameworks with $|\hat T|=|\hat S|$ and such that each $(\hat G,\hat t_i)\in \hat T$  is in canonical position with respect to $v_1,x,y$, and satisfies $\|\hat t_i-\hat q_i\| < \epsilon$ for  all $1\leq i\leq |\hat T|$.

Let $T$ be the set of all frameworks $(G,t_i)$ where $t_i=\hat t_i|_G$ and $(\hat G,\hat t_i)\in \hat T$.  We will show that, for sufficiently small $\epsilon$, we have $|T|=|\hat T|$ and each $(G,t_i)\in T$ is infinitesimally rigid and equivalent to $(G,t_1)$. Recall that $\hat q_i(v_2)=(0,0)$ and $\hat q_i(w)=(\pm a,0)$ for all $1\leq i\leq |\hat S|$. Let $\hat t_1(v_2)=(a_2,b_2)$ and $\hat t_1(w)=(a',b')$. Then, for sufficiently small $\epsilon$,  the fact that $\{v_1,x,w\}$ induce a triangle in $\hat G$ implies that $\hat t_i(w)=(\pm a',b')$. The fact that $\{v_2,x,w\}$ induces a triangle in $\hat G$ now implies that $\hat t_i(v_2)=(a_2,b_2)$ when $\hat t_i(w)=(a',b')$ and $\hat t_i(v_2)=(-a_2,b_2)$ when $\hat t_i(w)=(-a',b')$. This gives
$$d(\hat t_i(v_2)-t_i(v_1))=a_2^2+b_2^2=d(\hat t_1(v_2)-t_1(v_1))$$
so each $(G,t_i)$ is equivalent to $(G,t_1)$. The assertion that each $(G,t_i)$ is infinitesimally rigid now follows from the facts that $G$ is rigid and  $(G,t_i)$ is quasi-generic.

It remains to show that $|T|=|\hat T|$. Choose $i,j$ with $1\leq i<j\leq |\hat T|$. Since $\hat q_i\neq \hat q_j$, we have $\hat q_i(u)\neq \hat q_j(u)$ for some vertex $u$ of $\hat G$. If $u\neq w$ then the fact that $\hat t_i(u)$ and $\hat t_j(u)$ can be chosen to be arbitrarily close to $\hat q_i(u)$ and $\hat q_j(u)$, respectively, means we can ensure that $t_i(u)=\hat t_i(u)\neq \hat t_j(u)=t_j(u)$. Hence suppose that $u=w$. Interchanging $i,j$ if necessary, we have $\hat q_i(w)=(a,0)$ and $q_j(w)=(-a,0)$. This implies that $\hat t_i(w)=(a',b')$ and $\hat t_j(w)=(-a',b)$ and hence that $t_i(v_2)=\hat t_i(v_2)=(a_2,b_2)$ and $t_j(v_2)=\hat t_j(v_2)=(-a_2,b_2)$. Hence $t_i\neq t_j$ for all $1\leq i<j\leq |\hat T|$.

We can now combine the above inequalities to deduce that
$$r(G,t_1)\geq |T|=|\hat T|=|\hat S|=2|S|=2r(H,p)$$
and the result follows since $(G,t_1)$ is quasi-generic.
\eproof

\begin{theorem}\label{thm:realplanar}
Every  planar isostatic graph $G=(V,E)$ has a generic realisation $(G,p)$ such that $r(G,p) \geq 2^{|V|-3}$.
\end{theorem}

\bproof
Every isostatic planar graph can be reduced to a 3-cycle by a sequence of edge contractions  in such a way that each intermediate graph is planar and isostatic by \cite{FJW,OP}. Since each edge contraction reduces $|V|$ by one the result follows by induction using Theorem \ref{rvsplit} and the fact that $K_3$ is globally rigid.
\eproof

\section{Globally rigid graphs  and  globally linked pairs of vertices}

We first use Corollary \ref{h2} to characterise
graphs $G$ with $c(G)=1$. Our characterization is the same as
that given in \cite{JJ} for {globally rigid graphs} in $\real^2$. (This result can be deduced immediately from the characterisation in \cite{JJ} and the
result of Gortler and Thurston mentioned in the Introduction that generic global rigidity in $\real^d$ and $\complex^d$ are equivalent. We give our proof since it is short and direct.)


\begin{theorem}\label{globrigid}   Let $G=(V,E)$ be a graph with at least four vertices.
Then $c(G)=1$ if and only if $G$ is $3$-connected and redundantly rigid.
\end{theorem}
\bproof Necessity was proved for real (and hence also for complex)
generic realisations in \cite{H}. We prove sufficiency by
induction on $|V|+|E|$. If $G$ has four vertices then $G=K_4$ and
$c(G)=1$ since $G$ is complete. Hence suppose that $|V|\geq 5$. If
$G-e$ is 3-connected and redundantly rigid for some $e\in E$, then
$c(G-e)=1$ by induction, and hence $c(G)=1$. Thus we may suppose
that $G-e$ is not both 3-connected and redundantly rigid. By
\cite[Theorem 6.1]{JJ} there exists a vertex $v_n\in V$ with
$N(v)=\{v_1,v_2,v_3\}$ such that $H=G-v_n+v_1v_2$ is 3-connected and
redundantly rigid. This implies in particular that $G-v_n$ is rigid.
Induction and Corollary \ref{h2} now give $c(G)\leq c(H)=1$. \eproof

Let $(G,p)$ be a complex realisation of a rigid graph $G=(V,E)$ and
$u,v\in V$. We say that $\{u,v\}$ is {\em globally linked} in
$(G,p)$ if every equivalent complex realisation $(G,q)$ of $G$ has
$d(p(u)-p(v))=d(q(u)-q(v))$. It can be seen that ${u,v}$ is globally
linked in $(G,p)$ if and only if $c(G,p)=c(G+e,p)$, where $e=uv$.
Theorem \ref{c(G)} now implies that the property of being globally
linked is a {\em generic property} i.e. if $\{u,v\}$ is globally
linked in some generic complex realisation of $G$ then $\{u,v\}$ is
globally linked in all such realisations. We say that {\em $\{u,v\}$
is globally linked in $G$} if $\{u,v\}$ is globally linked in
some, or equivalently all, generic complex realisations of $G$.

The analogous concept for real realisations was introduced in
\cite{JJS}. (The situation for generic real realisations is more
complicated as it is not necessarily true that if $\{u,v\}$ is
globally linked in some generic real realisation of $G$ then
$\{u,v\}$ is globally linked in all generic real realisations.
For example the pair ${u,v}$ is globally linked in the real realisation in Figure \ref{fig3}, but
not in Figure \ref{fig4}.
This problem is circumvented in \cite{JJS} by defining {\em $\{u,v\}$ to be globally
linked in $G$ in $\real^2$} if $\{u,v\}$ is globally linked in {\em all}
generic real realisations of $G$.)

Our next result is analogous to a result for real realisations
given in \cite[Theorem 4.2]{JJS}.

\begin{theorem}\label{linkedextension}
Let $(G,p)$ be a generic complex realisation of a graph $G=(V,E)$
and $u,v,v_1,v_2,v_3,v_n \in V$ with $N(v_n)=\{v_1,v_2,v_3\}$ and
$v_n\neq u,v$. Let $H=G-v_n+v_1v_2$. Suppose that $G-v_n$ is rigid
and that $\{u,v\}$ is globally linked in $(H,p|_{H})$. Then
$\{u,v\}$ is globally linked in $(G,p)$.
\end{theorem}
\bproof Suppose $(G,q)$ is equivalent to $(G,p)$. Let
$p^*=p|_{H}$ and $q^*=q|_{H}$. Since $G-v_n=H-v_1v_2$ is
rigid, Lemma \ref{con3} implies that
$d(p(v_1)-p(v_2))=d(q(v_1)-q(v_2))$. Hence $(H,p^*)$ and $(H,q^*)$
are equivalent. Since $\{u,v\}$ is globally linked in $(H,p^*)$, we
have
$$d(p(u)-p(v))=d(p^*(u)-p^*(v))=d(q^*(u)-q^*(v))=d(q(u)-q(v)).$$
Thus $\{u,v\}$ is globally linked in $(G,p)$. \eproof

The real analogue of Theorem \ref{linkedextension} was used in
\cite[Section 5]{JJS} to characterize when two vertices in a generic
real realisation of an `$\scrrm$-connected graph' are globally linked
in $\real^2$. We can show that the same characterization holds for
complex realisations. We first need to introduce some new terminology.

A {\em matroid} ${\cal M}=(E,{\cal I})$, consists of a set $E$
together with a family $\cal I$ of subsets of $E$, called {\em
independent sets}, which satisfy three simple axioms which capture
the properties of linear independence in vector spaces, see
\cite{Oxley}. Given a complex realisation $(G,p)$ of a graph
$G=(V,E)$, its {\em rigidity matroid} ${\cal R}(G,p)=(E,{\cal I})$
is defined by taking $\cal I$ to be the family of all subsets of $E$
which correspond to linearly independent sets of rows in the
rigidity matrix of $(G,p)$. It is not difficult  to see that the set
of independent subsets of $E$ is the same for all generic complex
realisations of $G$. We refer to the resulting matroid as the {\em
rigidity matroid of $G$} and denote it by ${\cal R}(G)$.

Given a {\em matroid} ${\cal M}=(E,{\cal I})$ we may
define an equivalence relation
on $E$ by saying that $e,f\in E$ are related if $e=f$ or if there is
a {\em circuit}, i.e. minimal dependent set, $C$ of $\cal M$ with
$e,f\in C$. The
equivalence classes are called the {\it components} of $\cal M$. If
$\cal M$ has at least two elements and only one component then $\cal
M$ is said to be {\it connected}. We say that a graph $G=(V,E)$ is
{\em $\scrrm$-connected} if its rigidity matroid ${\cal R}(G)$ is
connected. The {\em $\scrrm$-components} of $G$ are the subgraphs of
$G$ induced by the components of ${\cal R}(G)$.
For more examples and basic properties of $\scrrm$-connected graphs
see \cite{JJ}. An efficient algorithm for constructing the
$\scrrm$-components of a graph is given in \cite{BJalgo}.

\begin{theorem}\label{linkedMcon}
Let $G=(V,E)$  be a  an $\scrrm$-connected graph and $u,v\in V$.
Then $\{u,v\}$ is globally linked in $G$ if and only if $u$ and $v$
are joined by three internally disjoint paths in $G$.
\end{theorem}
\bproof Necessity follows for real (and hence also complex)
generic realisations by \cite[Lemma 5.6]{JJS}.
Sufficiency
follows by applying the same proof technique as for \cite[Theorem
5.7]{JJS} but using Theorem \ref{linkedextension} in place of
\cite[Theorem 4.2]{JJS} \eproof

The following conjecture is a complex version of \cite[Conjecture
5.9]{JJS}. It would characterise when two vertices in a rigid graph
are globally linked.

\begin{con}\label{linkedcon}
Let $G=(V,E)$  be a rigid graph and $u,v\in V$. Then $\{u,v\}$ is
globally linked in $G$ if and only if either $uv\in E$ or $u$ and
$v$ are joined by three internally disjoint paths in some
$\scrrm$-connected component of $G$.
\end{con}
Note that the `sufficiency part' of Conjecture \ref{linkedcon}
follows from Theorem \ref{linkedMcon}.

\section{Separable graphs}\label{sep}

A $k$-{\em separation} of a graph $G=(V,E)$ is a pair $(G_1,G_2)$ of
edge-disjoint subgraphs of $G$ each with at least $k+1$ vertices
such that $G=G_1\cup G_2$ and $|V(G_1)\cap V(G_2)|=k$. If
$(G_1,G_2)$ is a $k$-separation of $G$, then we say that {\em $G$ is
$k$-separable} and that $V(G_1)\cap V(G_2)$ is a {\em $k$-separator}
of $G$. We will obtain expressions for $c(G)$ when $G$ is a rigid graph with
a 2-separation, and also when $G$ has a 3-separation induced by a 3-edge-cut.

\begin{lem}\label{cleavage}
Let $(G_1,G_2)$ be a $2$-separation of a rigid graph $G$ with
$V(G_1)\cap V(G_2)=\{v_1,v_2\}$ and let $H_i=G_i+e$ where $e=v_1v_2$
for $i=1,2$. Suppose that $\{v_1,v_2\}$ is globally linked in $G$.
Then $c(G)=2c(H_1)c(H_2)$.
\end{lem}
\bproof Let $(G,p)$ be a generic realisation of $G$ and choose
$d_0\in \complex$ with $d(p(v_1)-p(v_2))=d_0^2$ and $\Arg d_0\in
(0,\pi]$. Let $S$ be the set of all realisations $(G,q)$ which are
equivalent to $(G,p)$ and satisfy $q(v_1)=(0,0)$ and
$q(v_2)=(0,d_0)$. Lemma \ref{qgen1.6} and the hypothesis that
$\{v_1,v_2\}$ is globally linked in $G$ imply that
$|S|=2c(G)$.

The hypothesis that $G$ is rigid implies that both $H_1$ and $H_2$ are rigid.
For $i=1,2$, let $S_i$ be the set of all realisations $(H_i,q_i)$
which are equivalent to $(H_i,p|_{H_i})$ and satisfy
$q_i(v_1)=(0,0)$ and $q_i(v_2)=(0,d_0)$. Lemma \ref{qgen1.6} and
the fact that $v_1v_2\in E(H_i)$ imply that $|S_i|=2c(H_i)$. It is
straightforward to check that the map $\theta:S\to S_1\times S_2$
defined by $\theta(G,q)=[(H_1,q|_{V(H_1}),(H_2,q|_{V(H_2})]$ is a
bijection. Hence $2c(G)=|S|=|S_1|\times|S_2|=4c(H_1)c(H_2)$. \eproof

We next show that we can apply Lemma \ref{cleavage} when
$G$ has a 2-separation
$(G_1,G_2)$ in which $G_1$ and $G_2$ are both rigid. We need one
more piece of matroid terminology. An {\em $\scrrm$-circuit} in a
graph $G$ is a subgraph $H$ such that $E(H)$ is a circuit in the
rigidity matroid of $G$.

\begin{lem}\label{cleavage1}
Let $(G_1,G_2)$ be a $2$-separation of a rigid graph $G$ with
$V(G_1)\cap V(G_2)=\{v_1,v_2\}$ and let $H_i=G_i+e$ where $e=v_1v_2$
for $i=1,2$. Suppose that $G_1$ and $G_2$ are both rigid. Then
$\{u,v\}$ is globally linked in $G$ and $c(G)= 2c(H_1)c(H_2)$.
\end{lem}
\bproof We first show that $\{v_1,v_2\}$ is globally linked in $G$.
This holds trivially if $e\in E(G)$ and hence we may suppose that
$e\not\in E(G)$. Since $G_i$ is rigid, $e_i$ is contained in an
$\scrrm$-circuit $C_i$ of $H_i$ for each $i=1,2$. Then
$C=(C_1-e)\cup (C_2-e)$ is an $\scrrm$-circuit of $G$ by \cite[Lemma
4.1]{BJ}. We may now use Theorem
\ref{linkedMcon} to deduce that $\{u,v\}$ is globally linked in
$C$. Since $C\subseteq G$, $\{u,v\}$ is globally linked in $G$. The
fact that $c(G)= 2c(H_1)c(H_2)$ now follows immediately from Lemma
\ref{cleavage}. \eproof

In order to obtain results for graphs with 2-separations  $(G_1,G_2)$ in which $G_1$ and $G_2$
are not both rigid, we need a result concerning the number of complex
realisations of a rigid graph satisfying given `distance'
constraints.

\begin{lem}\label{realise} Let $G=(V,E)$ be a rigid graph
with $V=\{v_1,v_2,\ldots,v_n\}$, $E=\{e_1,e_2,\ldots,e_m\}$ and
$e_i=v_{i_1}v_{i_2}$ for all $1\leq i\leq m$. Suppose that
$T=\{e_1,e_2,\ldots,e_t\}\subseteq E$ is such that $\rankg(G-T)=\rankg(G)-t$.
Let $(G,p)$ be a generic realisation of $G$ and
$d^*_T=\{d^*_1,d^*_2,\ldots,d^*_t\}\subset \complex$ be algebraically
independent over $\rat(d_{G-T}(p))$. Then the number of pairwise
non-congruent realisations $(G,q)$ of $G$ with $(G-T,q)$ equivalent
to $(G-T,p)$ and $d(p(v_{i_1})-p(v_{i_2}))=d^*_i$ for all $e_i\in T$
is $c(G)$.
\end{lem}
\bproof Let $K=\rat(d_{G-T}(p))$. We will define polynomials $f_i\in
K[X,Y,D]$  for $1\leq i\leq m$, where $X=(X_1,X_2,\ldots,X_n)$, $Y=(Y_1,Y_2,\ldots,Y_n)$, and
$D=(D_1,D_2,\ldots,D_t)$ are indeterminates. We first associate two
variables $X_i,Y_i$ with each $v_i\in V$ and a variable $D_i$ with
each $e_i\in T$. We then put
$f_i=(X_{i_1}-X_{i_2})^2+(Y_{i_1}-Y_{i_2})^2-D_i$ for each $e_i\in
T$ and
$f_i=(X_{i_1}-X_{i_2})^2+(Y_{i_1}-Y_{i_2})^2-d(p(v_{i_1})-p(v_{i_2}))$
for each $e_i\in E\sm T$.

We now apply Lemma \ref{realise0}. We need to find $x,y\in
\complex^n$ and $d\in \complex^t$ such that $f_i(x,y,d)=0$ for all
$1\leq i \leq m$, and $\tran[K(d),K]=t$. This is easy since we can just put
$(x_i,y_i)=p(v_i)$ for all $v_i\in V$ and $d_i=d(p(v_{i_1})-p(v_{i_2}))$ for all $e_i\in T$, and
use the definition of the polynomials $f_i$ to deduce that $f_i(x,y,d)=0$
for all $1\leq i \leq m$. Since $G$ is rigid
$\tran[\rat(d_G(p)),\rat]=2n-3$ and
$\tran[\rat(d_{G-T}(p)),\rat]=\rankg(G-T)=2n-3-t$ by Lemma
\ref{grank}. Since
$\tran[\rat(d_G(p)),\rat]=\tran[K(d),K]+\tran[K,\rat]$ we have
$\tran[K(d),K]=t$. Since we also have $\tran[K(d^*_T),K]=t$, Lemma
\ref{realise0} implies that there exists a realisation $(G,q)$ with
$(G-T,q)$ equivalent to $(G-T,p)$ and $d(p(v_{i_1})-p(v_{i_2}))=d_i^*$
for all $e_i\in T$.

We may assume that $(G,q)$ is in canonical position with respect to $v_1,v_2,v_3$. Since $\rat(d_T^*)\subseteq \rat(d_G(q))$,
$\tran[\rat(d_G(q)):\rat]=\tran[\rat(d_G(q)):K]+\tran[K:\rat]
\geq
\tran[\rat(d_T^*):K]+\tran[K:\rat]\geq |T|+ 2n-3 -|T|=2n-3.$
 Since
$(G,q)$ is in canonical position and $\rat(d_G(q))\subseteq \rat(q)$
we must have $\tran[\rat(q):\rat]=2n-3$. We may now use Lemma
\ref{basic}(a) and (b) to construct a generic framework which is congruent to
$(G,q)$. Hence $(G,q)$ is quasi-generic and and the number of pairwise non-congruent realisations of $G$ which
are equivalent to $(G,q)$ is $c(G)$. \eproof

Our next result is needed to enable us to apply Lemma \ref{realise} to $k$-separations.

\begin{lem}\label{algind}
Let $H_1,H_2$ be rigid graphs. Put $H=H_1\cup H_2$, $H_3=H_1\cap
H_2$, and $T=E(H_3)$. Suppose that $H_3$ is isostatic and that
$\rankg(H_2-T)=\rankg(H_2)-|T|$.  Let $(H,p)$ be a quasi-generic
realisation of $H$, $G_1$ be a spanning rigid subgraph of $H_1$, and
$(G_1,q_1)$ be a realisation of $G_1$ which is equivalent to
$(G_1,p|_{G_1})$. Then $d^*_T=\{d(q_1(u)-q_1(v))\,:\,uv\in T\}$ is
algebraically independent over $\rat(d_{H_2-T}(p|_{H_2}))$.
\end{lem}
\bproof  If $T=\emptyset$ there is nothing to prove so we may
suppose that $|T|\geq 1$ and hence $|V(H_3)|\geq 2$. We may also
assume that $(H,p)$ and $(H_1,q_1)$ are both in canonical position
with $p(u)=(0,0)=q_1(u)$, $p(v)=(0,y)$ and $q_1(v)=(0,z)$ for some
$y,z\in \complex$ and some $u,v\in V(H_3)$.

Since $H_1,H_2$ are rigid, $H=H_1\cup H_2$ is rigid. Let $F$ be a spanning isostatic subgraph of $H$ which contains $T$
and let $F_i=F\cap H_i$. Then
\begin{eqnarray*}
|E(F)|&=&|E(F_1)|+|E(F_2)|-|T|\\
&\leq& (2|V(H_1)|-3)+2(|V(H_2)|-3)-2(|V(H_3)|-3)\\
&=&2|V(H)|-3.
\end{eqnarray*}
Equality must occur throughout and hence $F_i$ is a spanning
isostatic subgraph of $H_i$ for $i=1,2$. Lemma \ref{new} now implies
that
\begin{equation}\label{e7.1}
\overline{\rat(d_{F_1}(q_1))}=\overline{\rat(q_1)}=\overline{\rat(d_{G_1}(q_1))}
\end{equation}
and
\begin{equation}\label{e7.2}
\overline{\rat(d_{F_1}(p))}=\overline{\rat(p)}=\overline{\rat(d_{H_1}(p)))}
=\overline{\rat(d_{G_1}(p)))}.
\end{equation}
Since $(G_1,q_1)$ and $(G_1,p|_{G_1})$ are equivalent $d_{G_1}(q_1)=d_{G_1}(p)$.
Equations (\ref{e7.1}) and (\ref{e7.2})  now give
$\overline{\rat(d_{F_1}(p))}=\overline{\rat(d_{F_1}(q_1))}$
and hence
\begin{eqnarray*}
\overline{\rat(d_{H}(p))}&=&\overline{\rat(d_{H_1}(p),d_{H_2-T}(p))}\\
&=&\overline{\rat(d_{F_1}(p),d_{H_2-T}(p))}\\
&=&\overline{\rat(d_{F_1}(q_1),d_{H_2-T}(p))}.
\end{eqnarray*}
Thus
$$\tran[\rat(d_{H}(p)):\rat]=
\tran[\rat(d_{H_2-T}(p):\rat]+\tran[\rat(d_{F_1}(q_1):\rat(d_{H_2-T}(p)].
$$
By Lemma \ref{grank},
$\tran[\rat(d_{H}(p)):\rat]=\rankg(H)=2|V(H)|-3$ and
$$\tran[\rat(d_{H_2-T}(p)):\rat]=\rankg(H_2-T)=2|V(H_2)|-3-|T|.$$
Thus
\begin{eqnarray*}
\tran[\rat(d_{F_1}(q_1):\rat(d_{H_2-T}(p)]&=&2|V(H)|-3-(2|V(H_2)|-3-|T|)\\
&=&2|V(F_1)|-3=|E(F_1)|.
\end{eqnarray*}
Hence $d_{F_1}(q_1)$ is algebraically independent over
$\rat(d_{H_2-T}(p)$. Since $T\subseteq E(F_1)$, $d^*_T$ is
also algebraically independent over $\rat(d_{H_2-T}(p)$.
\eproof

\begin{lem}\label{cleavage2}
Let $(G_1,G_2)$ be a $2$-separation of a rigid graph $G$ with
$V(G_1)\cap V(G_2)=\{v_1,v_2\}$.  Suppose that $G_2$ is not rigid
and put $H_2=G_2+e$ where $e=v_1v_2$. Then $G_1$ and $H_2$ are both
rigid and $c(G)=2c(G_1)c(H_2)$.
\end{lem}
\bproof Let $F$ be a spanning isostatic subgraph of $G$. We have
$|E(F)\cap E(G_1)|\leq 2|V(G_1)|-3$, and $|E(F)\cap E(G_2)|\leq
2|V(G_2)|-4$ since $G_2$ is not rigid. Thus
\begin{eqnarray*}
|E(F)|&=&|E(F)\cap E(G_1)|+|E(F)\cap E(G_2)|\\
&\leq& 2|V(G_1)|-3+2|V(G_2)|-4 =2|V(F)|-3.
\end{eqnarray*}
Since $F$ is rigid, we must have
equality throughout. In particular $|E(F)\cap E(G_1)|= 2|V(G_1)|-3$
so $G_1$ is rigid.

Consider the $2$-separation $(G_1,H_2)$ of $H=G+e$, and let $F'$ be
a spanning isostatic subgraph of $H$ which contains $e$. Then
$|E(F')\cap E(H_2)|\leq 2|V(H_2)|-3$ and,  since $e\in E(F')$,
$|E(F')\cap E(G_1)|\leq 2|V(G_1)|-4$. Thus
\begin{eqnarray*}
|E(F')|&=&|E(F')\cap E(G_1)|+|E(F')\cap E(H_2)|\\
&\leq& 2|V(G_1)|-4+2|V(H_2)|-3=2|V(F')|-3.
\end{eqnarray*}
   Since $F'$ is
rigid, we must have equality throughout. In particular $|E(F')\cap
E(H_2)|= 2|V(H_2)|-3$ so $H_2$ is rigid.

Let $(G,p)$ be a generic realisation of $G$. For each $z\in
\complex\sm\{0\}$ with $\Arg z\in (0,\pi]$ let $S(z)$
   be the set of all realisations $(G,q)$ of
$G$ such that $(G,q)$ is equivalent to $(G,p)$, $q(v_1)=(0,0)$ and
$q(v_2)=(0,z)$. Define $S_1(z)$ and $S_2(z)$ similarly by replacing
$(G,p)$ by $(G_1,p|_{G_1})$ and $(H_2,p|_{H_2})$ respectively.
Lemma \ref{qgen1.6} and Theorem \ref{c(G)} imply that $S(z)$,
$S_1(z)$ and $S_2(z)$ are finite, and are non-empty for only
finitely many values of $z$. In addition we have
\begin{equation}\label{ec1.2}
2c(G)=\sum_{S(z)\neq \emptyset}|S(z)|\:\mbox{ and }\:
2c(G_1)=\sum_{S_1(z)\neq \emptyset}|S_1(z)|.\end{equation}

We will show that
\begin{equation}\label{ec1.3}
|S(z)|=2|S_1(z)|\, c(H_2)
\end{equation}
   for all
$z\in \complex\sm\{0\}$ with $\Arg z\in (0,\pi]$. If $S_1(z)=
\emptyset$ then we must also have $S(z)= \emptyset$, since for any
$(G,q)\in S(z)$ we would have $(G_1,q|_{V(G_1}))\in S_1(z)$, so
(\ref{ec1.3}) holds trivially.

We next consider the case when $S_1(z)\neq \emptyset$. Choose
$(G_1,q_1)\in S_1(z)$.  We may apply Lemma \ref{algind} with
$H=G+e$, $H_1=G_1+e$, $T=\{e\}$ and
$d^*_T=\{d(q_1(v_1)-q_1(v_2))\}$ to deduce that $d^*_T$ is algebraically
independent over $\rat(d_{H_2}(p))$. We may then apply Lemma
\ref{realise} (with $G=H_2$) and Lemma \ref{qgen1.6} to deduce that $|S_2(z)|=2c(H_2)$.
Since the map $\theta:S(z)\to S_1(z)\times S_2(z)$ by
$\theta(G,q)=[(G_1,q|_{V(G_1)}),(H_2,q|_{V(G_2)})]$ is a bijection,
we have
$$|S(z)|= |S_1(z)|\, |S_2(z)|= 2\,|S_1(z)|\, c(H_2).$$
Thus (\ref{ec1.3}) also holds when $S_1(z)\neq \emptyset$.

Equation (\ref{ec1.3}) and the fact that $c(H_2)\neq 0$ imply that
$S_1(z)= \emptyset$ if and only if $S(z)= \emptyset$. We can now use
equations (\ref{ec1.2}) and (\ref{ec1.3}) to deduce that

$$
c(G)=\sum_{S(z)\neq \emptyset}|S(z)|
=2\sum_{S_1(z)\neq \emptyset}|S_1(z)|\,c(H_2) = 2\,c(G_1)\,c(H_2).
$$
\eproof
Note that Lemma \ref{h1} is the special case of Lemma
\ref{cleavage2} when $G_2$ is a path of length two.

Lemmas \ref{cleavage} and \ref{cleavage2} immediately give

\begin{theorem}\label{thm:cleavage}
Suppose that $G$ is a rigid graph and $(G_1,G_2)$ is a $2$-separation of $G$ with $V(G_1)\cap V(G_2)=\{u,v\}$.
Then $G_1+uv$, $G_2+uv$, and at least one of $G_1,G_2$ are rigid. Furthermore:\\
(a) if $G_1$ and $G_2$ are both rigid then $c(G)=2c(G_1+uv)\,c(G_2+uv)$;\\
(b) if $G_1$ is rigid and $G_2$ is not rigid then $c(G)=2c(G_1)\,c(G_2+uv)$.
\end{theorem}

We next state a complementary result for $k$-separations when $k\geq 3$ and the common intersection is globally rigid. Its proof is straightforward.

\begin{theorem}\label{thm:cleavage_glob}
Suppose that $G$ is a rigid graph and $(G_1,G_2)$ is a $k$-separation of $G$ such that $k\geq 3$ and $G_1\cap G_2$ is globally rigid.
Then  $c(G)=2c(G_1)\,c(G_2)$.
\end{theorem}

\begin{figure}
\vspace{-2cm}
\centering
\hspace{-6cm}
\includegraphics[scale=0.7]{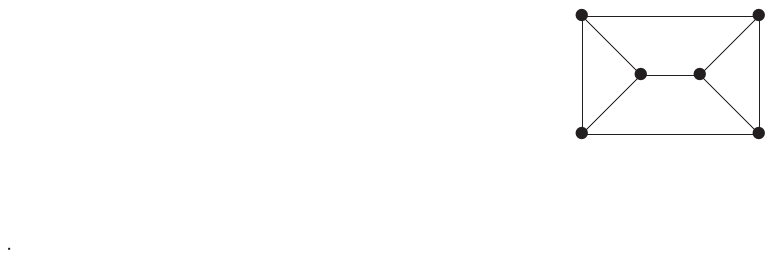}
\vspace{-15cm}
\caption{The triangular prism.} \label{fig0}
\end{figure}

We close this section by deriving a reduction formula for $c(G)$ when $G$ has a 3-edge-cut.
We first need to determine
$c(G)$ when $G$ is the triangular
prism
i.e. the
graph on six vertices consisting of two disjoint triangles joined by
a perfect matching shown in Figure \ref{fig0}.

\begin{lem}\label{prism}
Let $P$ be the triangular prism. Then $c(P)=12$.
\end{lem}
\bproof It is well known that every realisation of $P$ in $\complex^2$ has at most 12 equivalent, non-congruent realisations and that there exists a (real) realisation $(G,p)$ with 12 equivalent, non-congruent (real) realisations in which the vertices are not collinear, see for example \cite{BS,EP}. We can now use Theorem \ref{iso_gen_lem} to deduce that $c(G)=12$.
\eproof

\begin{theorem}\label{cleavage3}
Suppose that $G$ is a rigid graph and $G= G_1\cup G_2 \cup
\{e_1,e_2,e_3\}$ where $V(G_1)\cap V(G_2)=\emptyset$, $e_i=u_iv_i$
for $1\leq i\leq 3$, $u_1,u_2,u_3$ are distinct vertices of $G_1$,
and $v_1,v_2,v_3$ are distinct vertices of $G_2$. Then $G_1$ and
$G_2$ are rigid and $c(G)=12\,c(G_1)\,c(G_2)$.
\end{theorem}
\bproof Let $F$ be a spanning isostatic subgraph of $G$. We have
$|E(F)\cap E(G_1)|\leq 2|V(G_1)|-3$ and $|E(F)\cap E(G_2)|\leq
2|V(G_2)|-3$. Thus
\begin{eqnarray*}
|E(F)|&\leq&|E(F)\cap E(G_1)|+|E(F)\cap E(G_2)|+3\\
&\leq& 2|V(G_1)|-3+2|V(G_2)|-3+3 =2|V(F)|-3.
\end{eqnarray*}
Since $F$ is rigid, we must have equality throughout. In particular
$|E(F)\cap E(G_i)|= 2|V(G_i)|-3$ so $G_i$ is rigid for $i=1,2$.

\begin{clm}\label{c1}
Let $H_2$ be obtained from $G_2$ by adding the vertices
$u_1,u_2,u_3$ and edges $u_1u_2,u_2u_3,u_3u_1,u_1v_1,u_2v_2,u_3v_3$.
Then $c(G)= c(G_1)\,c(H_2)$.
\end{clm}
\bproof Let $(G,p)$ be a generic realisation of $G$. For each fixed
$b_2,a_3,b_3\in \complex\sm\{0\}$ with $\Arg b_2, \Arg a_3\in
(0,\pi]$ let $S(b_2,a_3,b_3)$
   be the set of all realisations $(G,q)$ of
$G$ such that $(G,q)$ is equivalent to $(G,p)$, $q(u_1)=(0,0)$,
$q(u_2)=(0,b_2)$ and $q(u_3)=(a_3,b_3)$. Define $S_1(b_2,a_3,b_3)$
and $S_2(b_2,a_3,b_3)$ similarly by replacing $(G,p)$ by
$(G_1,p|_{G_1})$ and $(H_2,p|_{H_2})$ respectively. Lemma
\ref{qgen1.6} and Theorem \ref{c(G)} imply that $S(b_2,a_3,b_3)$,
$S_1(b_2,a_3,b_3)$ and $S_2(b_2,a_3,b_3)$ are finite, and are
non-empty for only finitely many values of $b_2,a_3,b_3$. In
addition we have
\begin{equation}\label{ec2}
c(G)=\sum_{S(b_2,a_3,b_3)\neq \emptyset}|S(b_2,a_3,b_3)|\:\mbox{ and }\:
c(G_1)=\sum_{S_1(b_2,a_3,b_3)\neq \emptyset}|S_1(b_2,a_3,b_3)|.
\end{equation}

We
will show that
\begin{equation}\label{ec3}
|S(b_2,a_3,b_3)|=|S_1(b_2,a_3,b_3)|\, c(H_2)
\end{equation}
   for all
$b_2,a_3,b_3\in \complex\sm\{0\}$ with $\Arg b_2, \Arg a_3\in
(0,\pi]$. If $S_1(b_2,a_3,b_3)= \emptyset$ then we must
also have $S(b_2,a_3,b_3)= \emptyset$, since for any $(G,q)\in
S(b_2,a_3,b_3)$ we would have $(G_1,q|_{V(G_1}))\in S_1(b_2,a_3,b_3)$,
so (\ref{ec3}) holds trivially.

We next consider the case when $S_1(b_2,a_3,b_3)\neq \emptyset$.
Choose $(G_1,q_1)\in S_1(b_2,a_3,b_3)$. Let
$T=\{u_1u_2,u_2u_3,u_3u_1\}$ and
$d^*_T=\{d(q_1(u_i)-q_1(u_j))\,:\,u_iu_j\in T\}$. We may apply Lemma
\ref{algind} with $(H,p)=(G\cup T,p)$ and $(H_1,q_1)=(G_1\cup T,
q_1)$ to deduce that $d^*_T$ is algebraically independent over
$\rat(d_{H_2-T}(p))$. We may then apply Lemma
\ref{realise} (with $G=H_2$) to deduce that
$|S_2(b_2,a_3,b_3)|=c(H_2)$. Since the map
$\theta:S(b_2,a_3,b_3)\to S_1(b_2,a_3,b_3)\times S_2(b_2,a_3,b_3)$
by $\theta(G,q)=[(G_1,q|_{V(G_1)}),(G_2,q|_{V(G_2)})]$ is a
bijection, we have
$$|S(b_2,a_3,b_3)|= |S_1(b_2,a_3,b_3)|\, |S_2(b_2,a_3,b_3)|= |S_1(b_2,a_3,b_3)|\, c(H_2).$$
Thus (\ref{ec3}) also holds when $S_1(b_2,a_3,b_3)\neq \emptyset$.

Equation (\ref{ec3}) and the fact that $c(H_2)\neq 0$ imply that
$S_1(b_2,a_3,b_3)= \emptyset$ if and only if $S(b_2,a_3,b_3)= \emptyset$.
We can now use equations (\ref{ec2}) and (\ref{ec3}) to deduce that


$$
c(G)=\sum_{S(b_2,a_3,b_3)\neq
\emptyset}|S(b_2,a_3,b_3)|
=\sum_{S_1(b_2,a_3,b_3)\neq
\emptyset}|S_1(b_2,a_3,b_3)|\,c(H_2)
= c(G_1)c(H_2).
$$
This completes the proof of Claim \ref{c1}.
\eproof

We may apply the argument of Claim \ref{c1} to $H_2$ to deduce that
$c(H_2)=c(G_2)\,c(P)$, where $P$ is the triangular prism. Claim \ref{c1} and the fact that $c(P)=12$ now give $c(G)=12\,c(G_1)\,c(G_2)$. \eproof

\section{Two families of graphs}
We use the results from the previous section to determine $c(G)$ for
two important families of rigid graphs.

\subsubsection*{Quadratically solvable graphs}


Let $G=(V,E)$ be an isostatic graph with $E=\{e_1,e_2,\ldots,e_m\}$
and $e_i=u_iv_i$ for $1\leq i\leq m$. Then $G$ is  {\em
quadratically solvable} if for all $d=(d_1,d_2,...,d_m) \in
\complex^m$ such that $\{d_1,d_2,d_3,\ldots,d_m\}$ is algebraically
independent over $\rat$, there exists a realisation $(G,p)$ of $G$
with $d(p(u_i)-p(v_i))=d_i$ for all $1\leq i\leq m$, in which
$\rat(p)$ is contained in a  quadratic extension of $\rat(d)$ i.e.
there exists a sequence of field extensions $K_1\subset
K_2\subset\ldots\subset K_m$ such that $K_1=\rat(d)$, $K_m=\rat(p)$
and $K_{i+1}=K_i(x)$ for some $x^2\in K_i$ for all $1\leq i<m$.
These graphs are important in the theory of equation solving in
Computer Aided Design, see for example \cite{Gao,Ow}.

We may recursively construct an infinite family $\scrqs$ of
quadratically solvable isostatic graphs as follows. We first put the complete
graph on three vertices $K_3$ in $\scrqs$. Then, for any two graphs
$G_1,G_2\in \scrqs$, any two vertices $u_1, v_1$ in $G_1$, and any
edge $e=u_2v_2$ of $G_2$, we construct a new graph $G$ by `gluing'
$G_1$ and $G_2-e$ together along $u_1=u_2$ and $v_1=v_2$, and add
$G$ to $\scrqs$. The second author conjectured in \cite{Ow} that an
isostatic graph $G$ is quadratically solvable if and only if it belongs
to $\scrqs$. This conjecture was subsequently verified for isostatic
planar graphs in \cite{OP}.

\begin{theorem}\label{cQS}
Suppose $G\in \scrqs$. Then $c(G)=2^{|V(G)|-3}$.
\end{theorem}
\bproof We use induction on $|V(G)|$. If $|V(G)|=3$ then $G=K_3$ and
$c(G)=1$. Hence we may assume that $|V(G)|>3$. It follows from the
recursive definition of $\scrqs$ that there exists a 2-separation
$(G_1,G_2)$ of $G$ with $V(G_1)\cap V(G_2)=\{u,v\}$ and such that
$G_1$ and $G_2+uv$ both belong to $\scrqs$. The theorem now follows
from Lemma \ref{cleavage2} and induction. \eproof

We can use this result and the fact that all $\scrqs$ graphs can be constructed using vertex splits to determine the maximum number of generic real realisations for any $\scrqs$ graph. We will need a result from  \cite{OP} that a graph $G$ is in $\scrqs$ if and only if $G$ has a decomposition $G=G_1 \cup G_2 \cup G_3$ where $G_1 \cap G_2 = u_3$,
$G_2 \cap G_3 = u_1$, $G_3 \cap G_1 = u_2$ and each $G_i$ is either $K_2$ or is in $\scrqs$. We will also need the following concept:
 an edge $e$ in a $\scrqs$ graph $G$ is {\em contractible} if either $G/e\in \scrqs$ or $G/e=K_2$.

\begin{lem} \label{lem:qscontract}
Suppose that $G=(V,E)\in \scrqs$ and $|V|\geq 3$. Then $G$ has at least two contractible edges.
\end{lem}
\bproof
Induction on $|V|$. If $|V|=3$ then $G=K_3$ and $G$ has three contractible edges. Hence we may suppose that $|V|\geq 4$. By \cite{OP},  $G$ has a decomposition $G=G_1 \cup G_2 \cup G_3$ where $G_1 \cap G_2 = u_3$,
$G_2 \cap G_3 = u_1$, $G_3 \cap G_1 = u_2$ and each $G_i$ is either $K_2$ or in $\scrqs$. Since $|V|\geq 4$ we may assume that $G_1\neq K_2$. By induction, $G_1$ has two contractible edges $e$ and $f$. If neither $e$ nor $f$ is equal to $u_2u_3$ then they are both contractible in $G$ by \cite{OP} (since we have
$G/e=G_1/e \cup G_2 \cup G_3$). So suppose $e=u_2u_3$ and $f$ is contractible in $G$. If $|V(G_2)| \geq 3$ or $|V(G_3)|\geq 3$ then we can find another contractible edge in $G_2$ or $G_3$. Otherwise $G_2=K_2=G_3$ and the edges $u_1u_3$ and $u_1u_2$ are both contractible in $G$ (since $u_2u_3\in E$ and hence $G/u_1u_2=G/u_2u_3=G_1$).
\eproof

\begin{theorem}\label{rQS}
Suppose $G=(V,E)\in \scrqs$. Then the maximum value of $r(G,p)$ over all generic real realisations $(G,p)$ of $G$ is $2^{|V|-3}$.
\end{theorem}
\bproof
We use induction on $|V|$. If $|V(G)|=3$ then $G=K_3$ and
and the theorem holds, so we may assume that $|V|>3$.
Then $G$ has a contractible edge $e$ by Lemma \ref{lem:qscontract}.
We can now use Theorem \ref{rvsplit} and induction to deduce that $G$ has a generic real realisation $(G,p)$ such that $r(G,p) \geq 2^{|V|-3}$. On the other hand Theorem \ref{cQS} shows that $r(G,q) \leq 2^{|V|-3}$ for all generic real realisations $(G,q)$.
\eproof

\subsubsection*{$\scrrm$-connected graphs}

We will determine $c(G)$ when $G$ is an $\scrrm$-connected graph.
  We need some new
terminology. For each $\{u,v\}\subset V$, let $w_G(u,v)$ denote the
number of connected components of $G-\{u,v\}$ and put
$b(G)=\sum_{\{u,v\}\subset V}(w_G(u,v)-1)$. Note that $w_G(u,v)-1=0$
if $\{u,v\}$ is not a 2-separator of $G$, so we can assume that the
summation in the definition of $b(G)$ is restricted to pairs
$\{u,v\}$ which are 2-separators of $G$.

\begin{theorem}\label{cforM}
Let $G$ be an $\scrrm$-connected graph. Then $c(G)=2^{b(G)}$.
\end{theorem}
\bproof We use induction on $b(G)$. Suppose $b(G)=0$. Then $G$ is
3-connected and, since $G$ is $\scrrm$-connected, it is also
redundantly rigid. Hence $c(G)=1$ by Theorem \ref{globrigid}. Thus
we may assume that $b(G)\geq 1$.

Choose vertices $u,v$ of $G$ with $w_G(u,v)\geq 2$ and let
$(G_1,G_2)$ be a 2-separation in $G$ with
%
$V(G_1)\cap V(G_2)=\{u,v\}$. Let $H_i=G_i+uv$ for $i=1,2$. By
\cite[Lemma 5.3(b)]{JJS}, $H_i$ is $\scrrm$-connected for $i=1,2$.
In addition, \cite[Lemma 3.6]{JJ} implies that every 2-separator
$\{u',v'\}$ of $G$ which is distinct from $\{u,v\}$ is a 2-separator
of $H_i$ for exactly one value of $i\in \{1,2\}$, and,
for this value of $i$, satisfies $w_G(u',v')=w_{H_i}(u',v')$. Since we also have
$w_G(u,v)=w_{H_1}(u,v)+w_{H_2}(u,v)$, we may deduce that
$b(G)=b(H_1)+b(H_2)-1$. Using induction and Lemma \ref{cleavage} we
have
   $$c(G)=2\,c(H_1)\, c(H_2)=2\times 2^{b(H_1)}\times 2^{b(H_2)}=2^{b(G)}.$$
\eproof

Our expression for $c(G)$ in Theorem \ref{cforM}  is identical to
that given for $r(G,p)$ in \cite[Theorem 8.2]{JJS} when $(G,p)$ is a
generic real realisation of $G$, and provides an explanation for the
fact that $r(G,p)$ is the same for all generic real realisations
$(G,p)$ of an $\scrrm$-connected graph $G$.

\section{Closing Remarks and Open Problems}
The obvious open problem is:

\begin{prob}\label{p1} Can $c(G)$ be determined efficiently for an arbitrary rigid graph
$G$?
\end{prob}

\begin{figure}
\vspace{-2cm}
\centering
\includegraphics[scale=0.8]{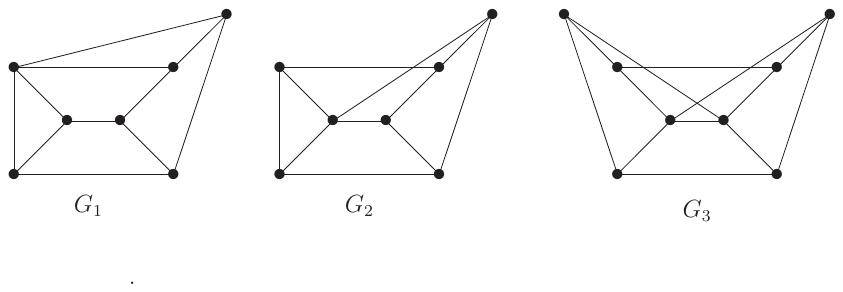}
\vspace{-18cm}
\caption{The graphs $G_1$, $G_2$ and $G_3$.} \label{fig1}
\end{figure}

Theorem \ref{cforM} gives an affirmative answer to this problem when
$G$ is $\scrrm$-connected and the results of Section \ref{sep} allow
us to reduce Problem \ref{p1} to the case when $G$ is 3-connected and all
3-edge-cuts of $G$ are trivial i.e. consist of three edges incident
to the same degree three vertex. On the other hand, the isostatic
graphs $G_1$, $G_2$ and $G_3$ of Figure \ref{fig1} indicate that it may be
difficult to obtain an affirmative answer to Problem \ref{p1} for
all graphs. Results from \cite{EP,EM} and Theorem \ref{iso_gen_lem} show that $c(G_1)=28$, and computer calculations indicate with high probability
that $c(G_2)=22$ and $c(G_3)=45$.  It is difficult to imagine how
these numbers could be deduced from the structures of $G_1$, $G_2$
and $G_3$.

If we cannot determine $c(G)$ precisely then we could ask for tight
asymptotic upper bounds on $c(G)$.

\begin{prob}\label{p2} Determine the smallest $k\in \real$ such that $c(G)=\mbox{O}(k^n)$ for all rigid
graphs $G$ with $n$ vertices.
\end{prob}

Clearly $c(G)$ will be maximised when $G$ is isostatic, and hence it
follows from \cite[Theorem 1.1]{BS} that
$c(G)\leq\frac{1}{2}{{2n-4}\choose{n-2}}\approx 4^n$ for all rigid
graphs $G$ with $n$ vertices. Borcea and Streinu \cite[Proposition
5.6]{BS} also construct an infinite family of isostatic graphs $G$
with $c(G)=12^{(n-3)/3}\approx 2.29^n$. Such a family can be
obtained by taking several copies of the triangular prism  with a single triangle in common.  The fact that
$c(G)=12^{(n-3)/3}$ for this family can be deduced from Lemmas \ref{thm:cleavage_glob} and \ref{triangle}.
Emiris and Moroz \cite{EM} use a similar construction using the graph $G_1$ to obtain an infinite family of isostatic graphs $G$
with $c(G)=28^{(n-3)/4}\approx 2.3^n$. It follows that the answer to Problem \ref{p2} will satisfy $28^{1/4}\leq k\leq 4$.

It would also be of interest to determine a tight lower bound on
$c(G)$ when $G$ is isostatic.

\begin{con}\label{iso} For all
isostatic graphs $G$ with $n$ vertices, $c(G)\geq 2^{n-3}$.
\end{con}
Note that this conjecture holds with equality for $\scrqs$ graphs  by Theorem \ref{cQS}, and also holds for planar graphs  by Theorem \ref{thm:realplanar}.

Since
every isostatic graph can be obtained from a triangle by type 1 or 2
Henneberg moves, and since every type 1 move doubles $c(G)$, it is tempting to try to prove Conjecture \ref{iso}
by showing that if we perform a type 2 move on an isostatic graph
$G$ then we will increase $c(G)$ by at least a factor of two.
Unfortunately this is (probably) not the case: the graph $G_2$ of
Figure \ref{fig1} can be obtained from the triangular prism $P$ by a
type 2 Henneberg move and we (probably) have $c(G_2)=22<2c(P)=24$.

We may also consider the problem of deciding which graphs have rigid realisations in $\complex^2$.

\begin{con}\label{iso_gen1} A graph $G=(V,E)$ has a realisation in $\complex^2$ which is rigid and has $d(p(u)-p(v))\neq 0$ for some $uv\in E$  if and only if $G$ is generically rigid.
\end{con}
It is straightforward to show that $G$ has a realisation in $\complex^2$ which is rigid (and has $d(p(u)-p(v))= 0$ for all $uv\in E$)  if and only if $G$ is connected.


Our final problem was posed by Dylan Thurston at a workshop on
global rigidity held at Cornell University in February 2011.

\begin{prob}\label{dylan} Does every rigid graph $G$ have a generic real
realisation $(G,p)$ such that $r(G,p)=c(G)$?
\end{prob}

The graph $G_3$ in Figure \ref{fig1} suggests that the answer to
this problem is most likely negative since the proof technique used
by Hendrickson \cite{H} to obtain necessary conditions for global
rigidity can be adapted to show that $r(G,p)$ is even for all
generic real realisations $(G,p)$ of a graph $G$ which is rigid but
not globally rigid.\footnote{Let $S$ be the set of all real realisations which are in
canonical position and are equivalent to $(G,p)$. If $G$ is not redundantly rigid then $G-e$
is not rigid for some edge $e$. In this case each component of the real configuration space of
$(G-e,p)$ will contain an even number of elements of $S$.  If $G$ is redundantly rigid then, since $G$
is not globally rigid, $G$ has a 2-separation. In this case reflecting one of the sides of the 2-separation
in the line through the two vertices of the corresponding 2-separator gives an involution on $S$ with no fixed points.}
On the other hand, we
(probably) have $c(G_3)=45$ which is odd.

We may say a bit more about this parity argument. Let $G=(V,E)$ be a
graph  which is rigid but not globally rigid and $S$ be the set of all realisations which are in
canonical position with respect to three given vertices $v_1,v_2,v_3$ and are equivalent to a given generic real
realisation $(G,p)$ of $G$. Since all edge lengths in $(G,p)$ are
real, the map $(G,q)\mapsto (G,q^*)$, where $q^*$ is obtained by
taking the complex conjugates of the coordinates of $q$ and then, if
necessary, reflecting the resulting framework in the axes to return
to canonical position, is an involution on $S$.

Suppose $(G,q^*)$ is equal to $(G,q)$ and let $q(v_1)=(0,0)$,
$q(v_2)=(0,y_2)$ and $q(v_3)=(x_3,y_3)$. Then $q^*(v_2)=(0,\pm
\bar{y}_2)=(0,y_2)$. Hence $y_2$ is either real or imaginary.

We
first consider the case when $y_2$ is real. We have $q^*(v_3)=(\pm
\bar{x_3},\bar{y}_3)=(x_3,y_3)$ so $x_3$ is either real or imaginary
and $y_3$ is real. If $x_3$ is real then we have $q^*(v_j)=(
\bar{x}_j,\bar{y}_j)=(x_j,y_j)$ for all $v_j\in V$ so $q$ is real.
If $x_3$ is imaginary then   $q^*(v_j)=(
-\bar{x}_j,\bar{y}_j)=(x_j,y_j)$ so $q(v_j)=(x_j,y_j)$ where $x_j$
is imaginary and $y_j$ is real for all $v_j\in V$.

We next consider
the the case when $y_2$ is imaginary. We have $q^*(v_3)=(\pm
\bar{x}_3,-\bar{y}_3)=(x_3,y_3)$ so $x_3$ is either real or
imaginary and $y_3$ is imaginary. If $x_3$ is imaginary then we have
$q^*(v_j)=( -\bar{x}_j,-\bar{y}_j)=(x_j,y_j)$ for all $v_j\in V$ so
$q$ is imaginary. This is impossible since $(G,q)$ is equivalent to
$(G,p)$ and so we must have $d(q(u)-q(v))>0$ for all $uv\in E$. If
$x_3$ is real then $q^*(v_j)=( \bar{x}_j,-\bar{y}_j)=(x_j,y_j)$ so
$q(v_j)=(x_j,y_j)$ where $x_j$ is real and $y_j$ is imaginary for
all $v_j\in V$.

In summary $(G,q^*)$ is equal to $(G,q)$ if and only if $q$ is
real, or we have $q(v_j)=(x_j,iy_j)$ where $x_j,y_j\in \real$ for
all $v_j\in V$, or we have $q(v_j)=(ix_j,y_j)$ where $x_j,y_j\in
\real$ for all $v_j\in V$. We will refer to the latter two such
realisations as {\em Minkowski realisations}.\footnote{We can
associate such realisations $q$ with  realisations $\tilde
q(v_j)=(x_j,y_j)$ in 2-dimensional Minkowski space where distance is
given by the Minkowski norm $d(x,y)=|-x^2+y^2|$. Results on generic global rigidity in $d$-dimensional Minkowski and other Pseudo-Euclidean spaces are given in \cite{GT}.} It follows that
the number of realisations in $S$ which are neither real nor
Minkowski must be even. As noted above, the number of real
realisations is also even. Thus it is the number of Minkowski
realisations which can be odd.

Although the answer to Problem \ref{dylan} seems to be negative, it
would still be of interest to find families of graphs for which the
answer is positive. For example Theorem \ref{cforM} and
\cite[Theorem 8.2]{JJS} give a positive answer when $G$ is
$\scrrm$-connected, and indeed show that $r(G,p)=c(G)$ for {\em all}
generic real realisations when $G$ is $\scrrm$-connected. Theorems
\ref{cQS} and \ref{rQS} show that we also have a positive answer
when $G\in \scrqs$.

\medskip
\noindent {\bf Acknowledgement} We would like to thank Shaun
Bullett, Peter Cameron,  and Bob Connelly for helpful conversations,
and the Fields Institute for support during its 2011 thematic
programme on Discrete Geometry and Applications.

\end{document}